\documentclass{tac}

\usepackage{amssymb,amsmath}
\usepackage{xy}
\xyoption{all}

\numberwithin{equation}{section}

\newcommand\V{\bigvee}
\newcommand\Max{\operatorname{Max}}
\newcommand\SMax[1]{\operatorname{{Max}^{#1}}}
\newcommand\oMax{\operatorname{Max}^{\circ}}
\newcommand\ie{i.e.}
\newcommand\eg{e.g.}
\newcommand\pwset[1]{\wp(#1)}
\newcommand\st{\mid}
\newcommand\cf{\textrm{cf.}}
\newcommand\diag{\operatorname{\triangle}}
\newcommand\specfunctor{\operatorname{\mathbf\Sigma}}
\newcommand\bunfunctor{\operatorname{\mathbf\Omega}}
\newcommand\rad{\operatorname{rad}}
\newcommand\spectrum{\operatorname{\Sigma}}
\newcommand\topology{\operatorname{\Omega}}
\newcommand\TopVect{\textit{TopVect}}
\newcommand\Vect{\textit{Vect}}
\newcommand\Loc{\textit{Loc}}
\newcommand\LLoc{\textit{LinLoc}}
\newcommand\sLLoc{\textit{sLinLoc}}
\newcommand\Top{\textit{Top}}
\newcommand\opp[1]{{#1}^{\textrm{op}}}
\newcommand\ident{\mathrm{id}}
\newcommand\CC{\mathbb{C}}
\newcommand\NN{\mathbb{N}}
\newcommand\SL{\mathit{SL}}
\newcommand\osupp{\operatorname{supp}^\circ}
\newcommand\Sub{\operatorname{Sub}}
\newcommand\oSub{\operatorname{Sub}^{\circ}}
\newcommand\spanmap{\operatorname{span}}
\newcommand\linspan[1]{\left\langle #1\right\rangle}
\newcommand\LBun{\textit{LinBun}}
\newcommand\QVBun{\textit{QVBun}}
\newcommand\sQVBun{\textit{sQVBun}}
\newcommand\pQVBun{\textit{QVBun}_{\spectrum}}
\newcommand\spQVBun{\textit{sQVBun}_{\spectrum}}
\newcommand\eval{\operatorname{eval}}

\newcommand\Gr{\operatorname{\boldsymbol{Gr}}}
\newcommand\GL{\operatorname{\boldsymbol{GL}}}
\newcommand\image{\operatorname{Im}}
\newcommand\interior{\operatorname{int}}
\newcommand\sob{\operatorname{sob}}
\newcommand\spat{\operatorname{spat}}
\newcommand\kfrak{\mathfrak k}
\newcommand\codim{\operatorname{codim}}
\newcommand\St{\boldsymbol V}
\newcommand\LM{\boldsymbol L}
\newcommand\upsegment{{\uparrow}}

\begin{document}

\title{Linear structures on locales}

\copyrightyear{2016}

\thanks{Work funded by FCT/Portugal through projects EXCL/MAT-GEO/0222/2012 and PEst-OE/EEI/LA0009/2013, and by COST (European Cooperation in Science and Technology) through COST Action MP1405 QSPACE.}

\author{Pedro Resende and Jo{\~a}o Paulo Santos}

\address{Centro de An\'alise Matem\'atica, Geometria e Sistemas Din\^amicos,
Departamento de Mate\-m\'{a}tica, Instituto Superior T\'{e}cnico,
Universidade de Lisboa,
Av.\ Rovisco Pais 1, 1049-001 Lisboa, Portugal}

\eaddress{pmr@math.tecnico.ulisboa.pt, jsantos@math.tecnico.ulisboa.pt}

\maketitle

\begin{abstract}
We define a notion of \emph{morphism} for quotient vector bundles that yields both a category $\QVBun$ and a contravariant global sections functor $C:\opp\QVBun\to\Vect$ whose restriction to trivial vector bundles with fiber $F$ coincides with the contravariant functor $\opp\Top\to\Vect$ of $F$-valued continuous functions. Based on this
we obtain a linear extension of the adjunction between the categories of topological spaces and locales: (i)~a linearized topological space is a \emph{spectral vector bundle}, by which is meant a mildly restricted type of quotient vector bundle; (ii)~a \emph{linearized locale} is a locale $\diag$ equipped with both a topological vector space $A$ and a $\diag$-valued \emph{support map} for the elements of $A$ satisfying a continuity condition relative to the spectrum of $\diag$ and the lower Vietoris topology on $\Sub A$; (iii)
we obtain an adjunction between the full subcategory of spectral vector bundles $\pQVBun$ and the category of linearized locales $\LLoc$, which restricts to an equivalence of categories between \emph{sober} spectral vector bundles and \emph{spatial} linearized locales. The spectral vector bundles are classified by a finer topology on $\Sub A$, called the \emph{open support topology}, but there is no notion of universal spectral vector bundle for an arbitrary topological vector space $A$.
\end{abstract}

\keywords{Quotient vector bundles, locales, Banach bundles, lower Vietoris topology, Fell topology}

\amsclass{06D22, 18B30, 18B99, 46A99, 46M20, 55R65}

\setcounter{tocdepth}{1}
\tableofcontents

\section{Introduction}

A quotient vector bundle \cite{RS} over a topological space $X$ is a triple $(\pi,A,q)$ that consists of a fiberwise linear quotient of a trivial complex vector bundle $\pi_2:A\times X\to X$ by an open map $q$, where $A$ can be any complex topological vector space:
\[
\xymatrix{
A\times X\ar[dr]_{\pi_2}\ar[rr]^q&&E\ar[dl]^{\pi}\\
&X
}
\]
The continuous and open map $\pi:E\to X$ is then a topological vector bundle of a very general kind, which in particular does not have to be locally trivial and can have variable or infinite rank.

Such a bundle has enough sections in the sense that for any $e\in E$ there is a continuous section $s:X\to E$ such that $e=s(\pi(e))$; for each $e=q(a,\pi(e))$ just take $s$ to be the section $\hat a$ defined by $\hat a(x)=q(a,x)$ for all $x\in X$. As an example, any Banach bundle $\pi:E\to X$ on a locally compact Hausdorff space $X$ can be made a quotient vector bundle by taking $A$ to be the space $C_0(\pi)$ of continuous sections that vanish at infinity, with the supremum norm topology, where $q:C_0(X)\times X\to E$ is the evaluation map $q(s,x)=\eval(s,x)=s(x)$. 

Quotient vector bundles are classified by continuous maps $\kappa:X\to\Sub A$, where $\Sub A$ is the set of all the linear subspaces of $A$ topologized with the lower Vietoris topology \cite{NT96,Vietoris} (with subclasses of bundles being classified by finer topologies such as the Fell topology). We shall refer to $\Sub A$ as the \emph{spectrum} of the topological vector space $A$. An aspect of $\Sub A$ has been ignored in \cite{RS}, namely its complete lattice structure, which has an important role to play. In order to introduce this idea, for each continuous section $s$ of $\pi$ let us write $\osupp s$ to denote the interior of the support of $s$,
\[
\osupp s = \interior\{x\in X\st s(x)\neq 0\}\;,
\]
and let us consider two maps,
\begin{eqnarray*}
\sigma&:&\Sub A\to\topology(X)\;,\\
\gamma&:&\topology(X)\to\Sub A\;,
\end{eqnarray*}
respectively called \emph{support map} and \emph{restriction map}, where $\topology(X)$ is the topology of $X$. These maps are
defined, for all $V\in\Sub A$ and $U\in\topology(X)$, by
\begin{eqnarray*}
\sigma (V) &=& \bigcup_{a\in V}\osupp \hat a\;;\\
\gamma (U)&=&\spanmap\{a\in A\st \osupp \hat a\subset U\}\;. 
\end{eqnarray*}
The topology $\topology(X)$ is another complete lattice under inclusion of open sets, and
it is clear that we have an equivalence
\[
\sigma (V)\subset U\iff V\subset \gamma (U)
\]
for all $V\in\Sub A$ and $U\in\topology(X)$. This means that the support map is left adjoint to the restriction map, and it immediately follows that they preserve suprema and infima, respectively:
\begin{eqnarray*}
\sigma\bigl(\spanmap\bigcup_\alpha V_\alpha\bigr)&=&\bigcup_\alpha\sigma (V_\alpha)\;,\\
\gamma\bigl({\operatorname{int}\bigcap_\alpha U_\alpha}\bigr)&=&\bigcap_\alpha \gamma(U_\alpha)\;.
\end{eqnarray*}

We shall regard the triple $(A,\sigma,\gamma)$ as a structure on the locale $\topology(X)$ in its own right, and, accordingly, we shall define the notion of \emph{linearized locale} $\mathfrak A=(\diag,A,\sigma,\gamma)$ to consist of a locale $\diag$ together with a topological vector space $A$ and a map $\sigma:\Sub A\to\diag$ with right adjoint  $\gamma$, such that the restriction of $\gamma$ to the prime spectrum $\spectrum(\diag)$ of $\diag$ is a continuous map into $\Sub A$. The latter means that the restricted map is the kernel map of a quotient vector bundle $\specfunctor(\mathfrak A)$, which we regard as a linear generalization of the spectrum of the locale $\diag$, and which, accordingly, we refer to as the \emph{spectrum} of $\mathfrak A$ (see \cite{BCL} for analogous terminology whereby a line bundle is the spectrum of a C*-category). The spectrum $\specfunctor(\mathfrak A)$ will be seen to be a quotient vector bundle of a special kind, in particular such that for all $a\in A$ the set
\[
\{p\in\spectrum(\diag)\st \hat a(p) \neq 0\}
\]
is open and thus coincides with the open support $\osupp\hat a$.

We shall see that from each quotient vector bundle $\mathcal A=(\pi:E\to X,A,q)$ we obtain a linearized locale $\bunfunctor(\mathcal A)=(\topology(X),A,\sigma,\gamma)$ as described above provided that the above openness condition holds with $\diag=\topology(X)$ for all $a\in A$, together with an additional continuity condition related to the spectrum of the locale $\topology(X)$: the restriction of $\gamma$ to the prime spectrum $\spectrum\topology(X)$ must be continuous. Hence, such quotient vector bundles will be termed \emph{spectral (quotient) vector bundles}.

We shall also compare spectral vector bundles and linearized locales by looking at $\specfunctor\bunfunctor(\mathcal A)$ and $\bunfunctor\specfunctor(\mathfrak A)$ for each spectral vector bundle $\mathcal A$ and each linearized locale $\mathfrak A$. In some cases we may obtain $\mathcal A\cong\specfunctor\bunfunctor(\mathcal A)$ and $\mathfrak A\cong\bunfunctor\specfunctor(\mathfrak A)$ according to obvious notions of isomorphism, but this leaves out many examples which, despite not yielding such isomorphisms, nevertheless ought to be comparable. Hence, in order to obtain a more satisfactory understanding of how quotient vector bundles and linearized locales relate, we shall define appropriate notions of morphism for both. In particular, our morphisms of bundles will differ from typical morphisms of vector bundles, which usually consist of pairs of maps $(f_0,f_1)$ that yield commutative squares
\[
\xymatrix{
F\ar[d]_{\rho}\ar[rr]^{f_1}&&E\ar[d]^{\pi}\\
Y\ar[rr]_{f_0}&&X
}
\]
but do not behave well with respect to global sections unless restrictions are imposed, such as local triviality of the bundles and requiring $f_1$ to be fiberwise a homeomorphism (as in \cite{Steenrod} for maps of fiber bundles). The morphisms in this paper subsume the latter and will still be such that a morphism $f:\rho\to\pi$ yields a contravariant linear map $f^*:C(\pi)\to C(\rho)$ on the spaces of global sections that in the case of trivial line bundles agrees with the contravariant Gelfand duality functor. This property enables us to obtain an adjunction between the category of spectral vector bundles and that of linearized locales, where $\bunfunctor$ is left adjoint to $\specfunctor$, which extends the classical adjunction between topological spaces and locales whereby the topology functor $\topology$ is left adjoint to the spectrum functor $\spectrum$ \cite{stonespaces}. Similarly to classical locale theory, we shall find corresponding notions of spatial linearized locale and of sober spectral vector bundles, and conclude that their respective categories are equivalent.

At the end of the paper we also show that the spectral vector bundles on suitable topological spaces $X$ (for instance sober spaces) are classified by the continuous maps $\kappa:X\to\oSub A$, where $\oSub A$ is the set $\Sub A$ with a topology herein called the \emph{open support topology}, which sits between the lower Vietoris and the Fell topologies. However, for other spaces $X$ the quotient vector bundles that correspond to continuous maps $\kappa:X\to\oSub A$ may fail to be spectral. In particular, in general there is no notion of universal spectral vector bundle for an arbitrary topological vector space $A$. In order to obtain some of these results we study the prime open sets of various topologies on $\Sub A$, in particular yielding a result of interest in its own right, namely that if $A$ is a locally convex space then $\Max A$, the space of closed linear subspaces of $A$ with the lower Vietoris topology, is a sober space.

Linearized locales are ``quotient vector bundles on locales''. We note that this is not a full-fledged localic notion, since it is formulated in terms of topological vector spaces. It is in principle possible to give a fully localic definition by considering localic vector spaces instead, and subsequently defining $\Sub A$ or $\Max A$ to be a localic sup-lattice as in \cite{RV}. However, we are not following such an approach, at least for now, because the motivation for our project stems from questions related to C*-algebras that we want to address without transforming them into questions about localic C*-algebras. Moreover, this choice will enable us to take advantage of quantale theory, without first having to develop a theory of localic quantales, when studying linearized locales $(\topology(G),A,\sigma,\gamma)$ such that $G$ is an open or \'etale groupoid and $A$ is a C*-algebra, in which case both $\Max A$ and $\topology(G)$ are quantales  \cite{MP1,MP2,KR,Re07,PR12}.

In a subsequent paper it will be seen that such linearized locales arise in connection with C*-algebras of groupoids \cite{Paterson}. In the case of an \'etale groupoid $G$ the space of units $G_0$ is open, and we obtain a pair $(A,B)$ consisting of the reduced C*-algebra $A=C_r^*(G)$ and a sub-C*-algebra $B=\gamma(G_0)$. More generally, this remains true for the reduced C*-algebra $C_r^*(E)$ obtained from sections of a Fell bundle $E$ on $G$ in the sense of \cite{Kumjian98}. For Fell line bundles on suitable groupoids these sub-C*-algebras are commutative and have been characterized precisely by Renault \cite{Renault} (see also \cite{Kumjian, RenaultLNMath,FeldmanMooreI-II}). A partial generalization to noncommutative subalgebras is due to Exel \cite{Exel} and it is formulated in terms of inverse semigroup Fell bundles \cite{Sieben}, whose relation to groupoid Fell bundles has been studied in \cite{BE12} in the semiabelian case, and, more recently, in \cite{BM16}.
The motivation for the present paper stems from the above program, in an attempt to provide general constructions via which the limitation to \'etale groupoids can be circumvented at least partially.

We thank the referee for his careful comments, which helped improve the presentation of this paper.

\section{Preliminaries on vector bundles}

We recall basic definitions from \cite{RS} in order to fix terminology and notation. 
Throughout this paper all the vector spaces are over $\CC$, and topological vector spaces satisfy no specific topological properties unless otherwise stated. Topological vector spaces satisfying the $T_1$ axiom are necessarily Hausdorff (in fact completely regular), and we refer to them as Hausdorff vector spaces. Such a space is finite dimensional if and only if it is locally compact, in which case it necessarily has the Euclidean topology. We denote by $\Vect$ (resp.\ $\TopVect$) the category of vector spaces (resp.\ topological vector spaces) and linear maps (resp.\ continuous linear maps). Given a subset $S\subset A$ of a vector space $A$ we denote the linear span of $S$ by $\spanmap S$, and if $S$ is a finite set $\{a_1,\ldots,a_k\}$ we also write $\linspan {a_1,\ldots,a_k}$ instead of $\spanmap\{a_1,\ldots,a_k\}$.

\subsection{Linear bundles}

\paragraph{Basic facts and definitions.}

Let $\pi:E\to X$ be a continuous map between topological spaces $E$ and $X$. For each $x\in X$ we refer to the set $\pi^{-1}(\{x\})$ as the \emph{fiber over $x$} and denote it by $E_x$.

By a \emph{linear structure} on $\pi$ will be meant a structure of vector space on each fiber $E_x$ such that the operations of scalar multiplication and vector addition are globally continuous when regarded as maps $\CC\times E\to E$ and $E\times_X E\to E$, respectively, and such that the \emph{zero section} of $\pi$, which sends each $x\in X$ to $0_x$ (the zero of $E_x$) is continuous.
Hence, $\pi$ equipped with a linear structure is a very loose form of vector bundle, which we refer to as a \emph{linear bundle}, the map $\pi$ itself being called its \emph{projection}.

The set $C(\pi)$ of continuous sections of $\pi$ is a vector space with the pointwise linear structure. The bundle is said to \emph{have enough sections} if the evaluation map
$\eval:C(\pi)\times X\to E$
is surjective.

If $A$ is a topological vector space, the projection $\pi_2:A\times X\to X$ will be called a \emph{trivial} linear bundle. Its space of continuous sections, $C(\pi_2)$, is linearly isomorphic to the vector space of continuous functions $C(X;A)$.

\paragraph{Pullbacks.}

Let $\pi:E\to X$ be a linear bundle, and let \[f:Y\to X\] be a continuous map. The pullback of $\pi$ along $f$ defines a linear bundle
\[
f^*(\pi) = f^*(E)\to Y\;,
\]
where $f^*(E)$ will usually be taken concretely to be the fibered product
\[
E\times_X Y = \{(e,y)\in E\times Y\st \pi(e)=f(y)\}\;,
\]
and $f^*(\pi)=\pi_2$. The pullback along $f$ induces contravariantly a linear map on continuous sections,
\[
s\mapsto f^*(s):C(\pi)\to C(\pi_2)\;,
\]
by, for all $s\in C(\pi)$ and $y\in Y$,
\[
f^*(s)(y)=(s(f(y)),y)\;.
\]

\paragraph{Morphisms over fixed base space.}

The category $\LBun(X)$ of linear bundles over $X$ has as \emph{morphisms} 
\[
(\rho:F\to X)\stackrel f\to(\pi:E\to X)
\]
the continuous and fiberwise linear maps $f:F\to E$ that commute with the projections:
\[
\xymatrix{
F\ar[dr]_{\rho}\ar[rr]^f&&E\ar[dl]^{\pi}\\
&X
}
\]
Such a morphism induces a linear map $s\mapsto f_*(s):C(\rho)\to C(\pi)$ by postcomposition: $f_*(s)=f\circ s$.

\subsection{Quotient vector bundles}

\paragraph{Basic definition and facts.}
By a \emph{quotient vector bundle} is meant a triple $\mathcal A = (\pi,A,q)$ where $\pi:E\to X$ is a linear bundle, $A$ is a topological vector space, and $q$ is a continuous open surjection that defines a morphism in $\LBun(X)$:
\[
\xymatrix{A\times X\ar[rr]^q\ar[rd]_{\pi_2}&&E\ar[dl]^{\pi}\\
&X}
\]
For each $x\in X$ we denote by $q_x:A\to E_x$ the continuous map defined for all $a\in A$ by
\[
q_x(a) = q(a,x)\;.
\]
This is a continuous open linear surjection, and thus every fiber $E_x$ is a quotient of $A$ as a topological vector space. Moreover, the quotient topology on $E_x$ coincides with the relative topology of $E_x$ as a subspace of $E$.

\paragraph{Sections.} For each $a\in A$ we have a continuous section of $\pi$,
\[
\hat a:X\to E\;,
\]
which is defined for all $x\in X$ by
\[
\hat a(x) = q(a,x)\;.
\]
Every $e\in E$ is of this form, and thus quotient vector bundles have enough sections.
The map $\widehat{(-)}:A\to C(\pi)$ thus defined is linear, and therefore its image $\widehat A$ with the quotient topology is a quotient topological vector space of $A$. The kernel of the quotient map is the \emph{radical of $\mathcal A$}:
\[
\rad\mathcal A = \bigcap_{x\in X}\ker q_x\;.
\]

\paragraph{Pullbacks.} If $Y$ is a topological space and $f:Y\to X$ is a continuous map, the \emph{pullback} $f^*(\mathcal A)$ is the quotient vector bundle $(f^*(\pi),A,f^*(q))$ such that for all $a\in A$ and $y\in Y$ the continuous open surjection $f^*(q):A\times Y\to E\times_X Y$ is defined by
\[
f^*(q)(a,y) = (q(a,f(y)),y)\;.
\]
Equivalently, $f^*(q)$ is the unique map that makes the following diagram commute:
\[
\xymatrix{
A\times Y\ar@{.>}[rd]|{f^*(q)}\ar[r]^-{\ident\times f}\ar@/_/[rdd]_{\pi_2}&A\times X\ar[dr]^q\\
&E\times_X Y\ar[d]^{f^*(\pi)}\ar[r]_-{\pi_1}&E\ar[d]^{\pi}\\
&Y\ar[r]_f&X
}
\]

\paragraph{Classifying spaces and universal bundles.} We denote by $\Sub A$ the collection of all the linear subspaces of $A$ topologized with the lower Vietoris topology \cite{NT96,Vietoris,RS}, which is generated by the sub basis whose open sets are, for each open set $U\subset A$,
\[
\widetilde U = \{P\in \Sub A \st P\cap U\neq\emptyset\}\;.
\]
We refer to $\Sub A$ as the \emph{spectrum} of $A$.
The \emph{kernel map} of $\mathcal A$ is the continuous map
\[
\kappa: X\to\Sub A
\]
which is defined for all $x\in X$ by $\kappa(x)=\ker q_x$.

Quotient vector bundles are classified by their kernel maps: from any continuous map
\[
\kappa:X\to\Sub A
\]
we obtain a quotient vector bundle $(\pi:E\to X,A,q)$ such that $E$ is the quotient of $A\times X$ by the equivalence relation $\sim$ which is defined for all $x,y\in X$ and $a,b\in A$ by
\[
(a,x)\sim (b,y)\iff (x=y\textrm{ and }a-b\in\kappa(x))\;.
\]
Then $q:A\times X\to E$ is the quotient map and $\pi$ is the factorization of $\pi_2:A\times X\to X$ through $q$. The kernel map of this bundle is $\kappa$.

If $\kappa:X\to\Sub A$ is the kernel map of a quotient vector bundle $(\pi:E\to X,A,q)$ the construction just described yields a quotient vector bundle $(\pi':E'\to X,A,q')$ such that $\pi$ is isomorphic to $\pi'$ in $\LBun(X)$ via an isomorphism $i:E\to E'$ that also commutes with $q$ and $q'$.

Finally we refer to the bundle obtained in this way from the identity map $\iota:\Sub A\to \Sub A$ as the \emph{universal quotient vector bundle for $A$}, and denote it by $UA=(\pi_A:EA\to \Sub A,A,q_A)$. Every quotient vector bundle for $A$ arises, up to an isomorphism as just described, as the pullback of $UA$ along the kernel map of the bundle.

\section{Categories of vector bundles}

In this section we study morphisms, over variable base space, of linear bundles and of quotient vector bundles, keeping in mind the general aim of obtaining a contravariant global sections functor.

\subsection{Contravariant morphisms of linear bundles}

\paragraph{Covariant morphisms.}

Let $\rho:F\to Y$ and $\pi:E\to X$ be linear bundles. There is more than one way in which to define a morphism from $\rho$ to $\pi$. An obvious choice is to let a morphism consist of a pair $(f_0,f_1)$ of continuous maps $f_0:Y\to X$ and $f_1:F\to E$ such that the following diagram commutes and $f_1$ is fiberwise linear:
\[
\xymatrix{
F\ar[rr]^{f_1}\ar[d]_{\rho}&&E\ar[d]^{\pi}\\
Y\ar[rr]_{f_0}&&X
}
\]
We call such a morphism \emph{covariant}, and composition is defined by $(g_0,g_1)\circ(f_0,f_1)=(g_0\circ f_0,g_1\circ f_1)$.

The sections of the bundles do not behave well with respect to such morphisms unless further conditions are imposed, as in the following example:

\example\label{coordinatebundles}
Let $\rho:F\to Y$ and $\pi:E\to X$ be locally trivial vector bundles with the same constant finite rank, and let $(f_0,f_1):\rho\to\pi$ be a covariant morphism such that $f_1$ is fiberwise a homeomorphism. Then a pullback map on sections $C(\pi)\to C(\rho)$ is obtained as in general for maps of coordinate bundles (\cf\ \cite[Lem.\ 2.11]{Steenrod}) by assigning to each $s\in C(\pi)$ the section $s'\in C(\rho)$ such that $s'(y) = f_y^{-1}(s(f_0(y)))$ for all $y\in Y$, where $f_y:E_y\to E_{f_0(y)}$ is the restriction of $f_1$ to $E_y$.
\endexample

In general a pushforward map does not exist:

\example
Let $Y=\CC^2\setminus\{(0,0)\}$ and let $\rho=\pi_1:F=Y\times\CC\to Y$ 
be the trivial line bundle over $Y$. Then $\CC\setminus\{0\}$ acts
on $F$ by multiplication, 
$\lambda\cdot(z_1,z_2,z_3)=(\lambda z_1,\lambda z_2,\lambda z_3)$,
and similarly on $Y$. Denoting the quotients by these actions
by
\begin{align*}
f_1&\colon F\to \CC P^2\setminus\{[0,0,1]\}\\
f_0&\colon Y\to \CC P^1
\end{align*}
we obtain a covariant morphism $(f_0,f_1)$ from $\rho$ to the Hopf line bundle \[\pi:\CC P^2\setminus\{[0,0,1]\}\to \CC P^1\;.\] The latter is nontrivial line bundle, so every 
section of $\pi$ has at least one zero. Since
$f_1$ is fiberwise a homeomorphism 
and $\rho$ has nowhere vanishing sections,
there can be no induced map
from the sections of $\rho$ to the sections of $\pi$.
\endexample

A way to understand the relation of covariant morphisms to sections is to note that $f_1$ determines a unique continuous map $f_1'$ that makes the following diagram commute:
\[
\xymatrix{
F\ar@/^/[rrrd]^{f_1}\ar@/_/[ddr]_{\rho}\ar@{.>}[rd]|{f_1'}\\
&E\times_X Y\ar[d]_{\pi_2}\ar[rr]_{\pi_1}&&E\ar[d]^{\pi}\\
&Y\ar[rr]_{f_0}&&X}
\]
The map $f_1'$ is a morphism in $\LBun(Y)$ and, conversely, any such $f_1'$ defines a map $f_1$ such that $(f_0,f_1)$ is a covariant morphism. Hence, $(f_0,f_1)$ could have been defined equivalently to be the pair $(f_0,f_1')$.
The map $f_1'$ yields a linear map $(f_1')_*:C(\rho)\to C(\pi_2)$ on sections, whereas the pullback gives us a linear map $f_0^*:C(\pi)\to C(\pi_2)$. So in general we do not have a map relating $C(\rho)$ and $C(\pi)$. But this discussion shows that a natural way of obtaining a pullback of sections along $f_0$ is, instead of placing strong constraints on the bundles and on $f_1$ as in Example \ref{coordinatebundles}, simply to reverse the direction of $f_1'$, as we see next.

\paragraph{Contravariant morphisms.}

Let $\rho:F\to Y$ and $\pi:E\to X$ be linear bundles. By a \emph{contravariant morphism} $f:\rho\to\pi$ is meant a pair of continuous maps $(f_\flat,f^\sharp)$ such that $f^\sharp:f_\flat^*(\pi)\to \rho$ is a morphism in $\LBun(Y)$, as follows:
\[
\xymatrix{
F\ar[d]_{\rho}&&E\times_X Y\ar[ll]_{f^\sharp}\ar[dll]_{\pi_2}\ar[rr]^{\pi_1}&&E\ar[d]^{\pi}\\
Y\ar[rrrr]_{f_\flat}&&&&X
}
\]
So fiberwise $f^\sharp$ induces linear maps $f^\sharp_y:E_{f_\flat(y)}\to F_y$ for each $y\in Y$.
\example\label{exampleT*M}
If $X$ and $Y$ are smooth manifolds and $E=T^*X$, $F=T^*Y$
are the (complexified) cotangent bundles, any smooth map $f_\flat\colon Y\to X$
induces a contravariant morphism $(f_\flat,f^\sharp)$ with
$f^\sharp=f_\flat^*$.
\endexample

In order to compose contravariant morphisms one should think of the pair $(f^\sharp,\pi_1)$ as a ``binary relation'' between $F$ and $E$ (a span), and composition should be defined accordingly as composition of spans by pullback:
\[
\xymatrix{
&&~\ar[dl]\ar[dr]\ar@/_30pt/[ddll]_{(f\circ g)^\sharp}\\
&~\ar[dl]_{g^\sharp}\ar[dr]^{\pi_1}&&~\ar[dl]_{f^\sharp}\ar[dr]^{\pi_1}\\
&&&&
}
\]
Although composition of spans is in general defined only up to isomorphism, there is a definition for contravariant morphisms which is associative ``on the nose''. Let $\sigma:G\to Z$, $\rho:F\to Y$, and $\pi:E\to X$ be linear bundles, with contravariant morphisms as follows:
\[
\xymatrix{
\sigma\ar[r]^g&\rho\ar[r]^f&\pi
}
\]
Then the \emph{composition} $f\circ g$ is the pair $(f_\flat\circ g_\flat, (f\circ g)^\sharp)$, where $(f\circ g)^\sharp:E\times_X Z\to G$ is defined, for all $e\in E$ and $z\in Z$, by
\begin{equation}\label{spancomposition}
(f\circ g)^\sharp(e,z) = g^\sharp(f^\sharp(e,g_\flat(z)),z)\;.
\end{equation}
We note that this definition depends on all the four continuous maps involved. This is explicit in \eqref{spancomposition} for $g_\flat$, $g^\sharp$, and $f^\sharp$, and the dependence on $f_\flat$ is via the definition of the pullback $E\times_X Z$:
\[
\xymatrix{
E\times_X Z\ar[rr]^{\pi_1}\ar[d]_{\pi_2}&&E\ar[d]^{\pi}\\
Z\ar[rr]_{f_\flat\circ g_\flat}&&X
}
\]
It is now straightforward to verify that the composition is associative and that the following \emph{identity morphism} on $\pi:E\to X$ is a neutral element for composition both on the right and on the left:
\[
\ident_\pi=(\ident_X,\pi_1:E\times_X X\stackrel{\cong}{\to} E)\;.
\]
We denote the resulting category by $\LBun$. A morphism $f$ is an isomorphism in this category if and only if both $f_\flat$ and $f^\sharp$ are homeomorphisms.
Note that there is always a contravariant morphism $(f,\ident):f^*(\pi)\to\pi$.

We remark that for each topological vector space $A$ there is a faithful functor $B_A:\Top\to\LBun$ that to each topological space $X$ assigns the trivial bundle
\[B_A(X)=\pi_2:A\times X\to X\]
and to each continuous map $f:Y\to X$ assigns the contravariant morphism
\[
B_A(f) = (f,(A\times X)\times_X Y\stackrel{\cong}{\to}{A\times Y})\;.
\]

\paragraph{Functoriality of sections.} Let $\rho:F\to Y$ and $\pi:E\to X$ be linear bundles, and $f:\rho\to\pi$ a contravariant morphism.
The composition
\[
\widehat f:=f^\sharp_*\circ f_\flat^*:C(\pi)\to C(\rho)
\]
is a linear map, and thus we have a contravariant functor
\[
C:\opp\LBun\to\Vect
\]
defined on linear bundles $\pi$ and contravariant morphisms $f$ by
\begin{eqnarray*}
\pi&\mapsto&C(\pi)\;,\\
f&\mapsto& \widehat f\;.
\end{eqnarray*}
This generalizes, for each topological vector space $A$, the contravariant functor
\[
C({-};A):\opp\Top\to\Vect
\]
in the sense that $C(-;A) = C\circ B_A$.

\example\label{coordinatebundles2}
Let $\rho:F\to Y$ and $\pi:E\to X$ be locally trivial vector bundles with the same constant finite rank, and let $(f_0,f_1):\rho\to\pi$ be a covariant morphism such that $f_1$ is fiberwise a homeomorphism. Then a contravariant morphism $f=(f_\flat,f^\sharp)$ is obtained by setting $f_\flat=f_0$ and $f^\sharp(e,y)= f_y^{-1}(e)$ for all $(e,y)\in E\times_X Y$, where $f_y:E_y\to E_{f_\flat(y)}$ is the restriction of $f_1$ to $E_y$, and the continuity of $f^\sharp$ follows from local triviality. The map $\widehat f$ is the same as the maps of sections in Example \ref{coordinatebundles}.
\endexample

We note that if $f:\rho\to\pi$ is a morphism in $\LBun$ and $\pi$ has enough sections then $f^\sharp$ is uniquely determined by $f_\flat$ and $\widehat f$, since for all $(e,y)\in E\times_X Y$ we have $e=s(f_\flat(y))$ for some $s\in C(\pi)$ and
\begin{equation}\label{sharpbystar}
f^\sharp(e,y)=f^\sharp(s(f_\flat(y)),y)=f^\sharp_*(f_\flat^*(s))(y)=\widehat f(s)(y)\;.
\end{equation}
A similar phenomenon occurs in sheaf theory, where the notion of cohomomorphism for sheaves of vector spaces is equivalent to that of contravariant morphism between \'etale linear bundles, as we show in the following example.

\example\label{cohomomorphisms}
Let the linear bundles $\rho:F\to Y$ and $\pi:E\to X$ be local homeomorphisms, and let $f:\rho\to\pi$ be a contravariant morphism. The map on sections $\widehat f$ applies equally to local sections, so for each continuous local section $s:U\to E$ of $\pi$ we obtain a continuous local section $\widehat f(s):f_\flat^{-1}(U)\to F$ of $\rho$. This property means that $f_\flat$ together with the induced maps on the fibers $f^\sharp_y:E_{f_\flat(y)}\to F_y$ is a cohomomorphism of sheaves \cite{Bredon}. Conversely, from any cohomomorphism $(f_\flat,(f^\sharp_y)_{y\in Y})$ with linear fiber maps $f^\sharp_y$ we obtain a contravariant morphism such that $f^\sharp$ is given by $f^\sharp(e,y)=f^\sharp_y(e)$ and is easily seen to be continuous due to the existence of enough continuous local sections and the fact that for a local homeomorphism the local sections are open maps: if $U\subset F$ is open we have, writing $s$ generically for a local section of $\pi$,
\begin{eqnarray*}
(f^\sharp)^{-1}(U) &=& \{(e,y)\in E\times_X Y\st f^\sharp_y(e)\in U\}\\
&=& \bigcup_s \{(s(f_\flat(y)),y)\st f^\sharp_y(s(f_\flat(y)))\in U\}\\
&=& \bigcup_s\{f^*(s)(y)\st \widehat f(s)(y)\in U\}\\
&=& \bigcup_s f^*(s)(\widehat f(s)^{-1}(U))\;,
\end{eqnarray*}
and this is open because $\widehat f(s)$ is continuous by hypothesis.
\endexample

\subsection{The category of quotient vector bundles}

\paragraph{Contravariant morphisms revisited.} Let us look again at contravariant morphisms, now for linear bundles which are equipped with the additional structure of a quotient vector bundle. The notion of morphism which we are about to define can be regarded, in our context, as the natural adaptation of the notion of cohomomorphism (\cf\ Example~\ref{cohomomorphisms}).

Let $\mathcal A=(\pi:E\to X, A, q)$ and $\mathcal B=(\rho:F\to Y, B, r)$ be quotient vector bundles. By a \emph{morphism} $f:\mathcal B\to\mathcal A$ will be meant a contravariant morphism of linear bundles $(f_\flat,f^\sharp):\rho\to\pi$ whose induced map on sections $\widehat f\colon C(\pi)\to C(\rho)$ satisfies
\begin{equation}\label{hatcont}
\widehat f(\widehat A)\subset\widehat B\;,
\end{equation}
together with a continuous linear map $f^*:A\to B$ such that the following diagram commutes:
\begin{equation}\label{morphismcommutation}
\vcenter{
\xymatrix{
B\ar@{->>}[d]_{\widehat{(-)}}&& A\ar@{->>}[d]^{\widehat{(-)}}\ar[ll]_{f^*}\\
\widehat B&&\widehat A\ar[ll]^{\widehat f}
}}
\end{equation}
In particular, this implies the following condition on radicals:
\begin{equation}\label{radicalcondition}
f^*(\rad\mathcal A)\subset \rad\mathcal B
\;.
\end{equation}

These morphisms form a category $\QVBun$, which we refer to as the \emph{category of quotient vector bundles}:
\begin{itemize}
\item The composition $f\circ g$ is the pair consisting of the contravariant morphism $(f_\flat,f^\sharp)\circ(g_\flat,g^\sharp)$ together with
the continuous linear map $g^*\circ f^*$;
\item The identity morphism on $\mathcal A=(\pi:A\to X,A,q)$ is the identity contravariant morphism on $\pi$ together with $\ident_A$.
\end{itemize}
We note that a morphism $f$ in this category is an isomorphism if and only if its three components are homeomorphisms.

\lemma\label{uniquefsharp}
Let $\mathcal A=(\pi:E\to X, A, q)$ and $\mathcal B=(\rho:F\to Y, B, r)$ be quotient vector bundles. If $f:\mathcal B\to\mathcal A$ is a morphism in $\QVBun$ then $f^\sharp$ is uniquely defined for all $a\in A$ and $y\in Y$ by
\[
f^\sharp(q(a,f_\flat(y))=r(f^*(a),y)\;.
\]
\endlemma

\proof
For all $(e,y)\in E\times_X Y$ we have $e=\hat a(f_\flat(y))$ for some $a\in A$ and thus, by \eqref{sharpbystar} and \eqref{morphismcommutation},
\[
f^\sharp(q(a,f_\flat(y))=f^\sharp(e,y)=\widehat f(\hat a)(y)=\widehat{f^*(a)}(y)=r(f^*(a),y)\;. 
\]
\endproof

There is a similar result
when the map $B\to\widehat B$ is a homeomorphism: then the map $f^*$ is completely determined by the contravariant morphism $(f_\flat,f^\sharp)$. In this case we
can characterize the morphisms of quotient vector bundles as the contravariant morphisms of linear bundles which satisfy \eqref{hatcont}.

\example\label{examplemorphQVBun}
As mentioned in section~1, from a Banach bundle 
over a locally compact Hausdorff space we can obtain
a quotient vector bundle by taking $A$ to be the space $C_0(\pi)$
of sections which vanish at infinity with the supremum norm topology.
In that case we have $A\cong\widehat A$.
Given a contravariant morphism of Banach bundles
$(f_\flat,f^\sharp)\colon \rho\to\pi$,
if $f_\flat$ is proper then $\widehat f\bigl(C_0(\pi)\bigr)\subset C_0(\rho)$,
so we obtain a morphism of the associated quotient vector bundles.
This applies in particular to Example~\ref{exampleT*M}
whenever $f_\flat$ is proper.
\endexample

\paragraph{Equivalent definition of morphism.}
From Lemma~\ref{uniquefsharp} it follows that if $f=(f_\flat,f^\sharp,f^*)$ is a morphism of quotient vector bundles then the following condition, relating $f_\flat$ and $f^*$, must hold for all $a\in A$ and $y\in Y$:
\begin{equation}\label{qvbmorphism}
q(a,f_\flat(y))=0\quad\Longrightarrow\quad r(f^*(a),y)=0\;.
\end{equation}
There is a converse to this:

\lemma
Let $\mathcal A=(\pi:E\to X, A, q)$ and $\mathcal B=(\rho:F\to Y, B, r)$ be quotient vector bundles. Every pair of maps
\begin{eqnarray*}
f_\flat&:&Y\to X\quad\textrm{(continuous)}\\
f^*&:&A\to B\quad\textrm{(continuous linear)}
\end{eqnarray*}
satisfying \eqref{qvbmorphism} arises from a morphism $f=(f_\flat,f^\sharp,f^*)$ in $\QVBun$ for a unique map $f^\sharp$.
\endlemma

\proof
Condition \eqref{qvbmorphism} implies that we may define a fiberwise linear map
\[
f^\sharp:E\times_X Y\to F
\]
for all $y\in Y$ and $e\in E_{f_\flat(y)}$ by
\begin{equation}\label{auxfsharp}
f^\sharp(e,y)=r(f^*(a),y)\;,
\end{equation}
where $a\in A$ is any element of $A$ such that $q(a,f_\flat(y))=e$.
In other words, $f^\sharp$ is the unique map that makes the following diagram commute:
\[
\xymatrix{
B\times Y\ar@{->>}[d]_{r}&&A\times Y\ar[ll]_-{f^*\times\ident}&(A\times X)\times_X Y\ar[l]_-{\cong}\ar@{->>}[d]^{q\times\ident}\\
F&&&E\times_X Y\ar@{.>}[lll]^-{f^\sharp}
}
\]
Now $q\times\ident$ is a quotient map because it is open and surjective. Therefore $f^\sharp$ is continuous, and thus we have a contravariant morphism $(f_\flat,f^\sharp):\rho\to\pi$.

Next we show that the triple $f=(f_\flat,f^\sharp,f^*)$ is a morphism in $\QVBun$.
Let $a\in A$ and $y\in Y$. We have $\widehat f = f^\sharp_*\circ f_\flat^*$, and thus, letting $e= q(a,f_\flat(y))$ in \eqref{auxfsharp}, we obtain
\[
\widehat{f^*(a)}(y) = r(f^*(a),y) = f^\sharp(q(a,f_\flat(y)),y) = f^\sharp(\hat a(f_\flat(y)),y) = \widehat f(\hat a)(y)\;,
\]
showing that \eqref{morphismcommutation} holds. 
\endproof

Hence, we see that the morphisms of quotient vector bundles can be defined alternatively but equivalently to be pairs $(f_\flat,f^*)$ satisfying \eqref{qvbmorphism}, with $f^\sharp$ being derived from the formula in Lemma~\ref{uniquefsharp}:
\[
f^\sharp(q(a,f_\flat(y))=r(f^*(a),y)\;.
\]
This formula further implies that $f^\sharp$ is fiberwise injective if and only if for all $a\in A$ and $y\in Y$ we have
\begin{equation}\label{qvbmorphism3}
q(a,f_\flat(y))=0\iff r(f^*(a),y)=0\;.
\end{equation}
A morphism of quotient vector bundles such that the above equivalence holds is called \emph{strict}, whereas if only the implication \eqref{qvbmorphism} holds it is called \emph{lax}. The strict morphisms yield a subcategory of quotient vector bundles, which we refer to as the \emph{strict category of quotient vector bundles}, denoted by $\sQVBun$.

\example
Let $\rho:F\to Y$ and $\pi:E\to X$ be locally trivial finite rank vector bundles over locally compact Hausdorff spaces, and let $(f_0,f_1)$ be a covariant morphism such that $f_0$ is a proper map and $f_1$ is fiberwise a homeomorphism. Then we obtain, by Example~\ref{coordinatebundles2} and Example~\ref{examplemorphQVBun}, a strict morphism of quotient vector bundles $(f_0,f^\sharp):(\rho,C_0(\rho),\eval)\to (\pi,C_0(\pi),\eval)$.
\endexample

\section{Preliminaries on pointfree topology}

In this section we provide some background on sup-lattices and locales. More details can be found for instance in \cite{JT,stonespaces}.

\subsection{Sup-lattices}

\paragraph{Basic definitions and examples.}

By a \emph{sup-lattice} \cite[Ch.\ I]{JT} will be meant a partially ordered set $L$ for which any subset $S$ has a join (supremum), which is denoted by $\V S$ or $\sup S$. Similarly, we write $\V_\alpha a_\alpha$ or $\sup_\alpha a_\alpha$ for the join of an indexed family $(a_\alpha)$ in $L$. Sup-lattices necessarily have meets (infima) of every subset (denoted by $\bigwedge S$ or $\inf S$), but the terminology we adopt reflects (as opposed to ``complete lattice'') the fact that joins are the first-class algebraic operations in the sense that they are preserved by \emph{homomorphisms} $f:L\to M$:
\[
f\left(\V S\right) = \V f(S)\;.
\]
The \emph{category of sup-lattices} $\SL$ is defined to have sup-lattices as objects and homomorphisms as arrows.
Additionally, we adopt the usual notation from lattice theory for meets and joins of pairs of elements, writing $a\vee b$ and $a\wedge b$ instead of $\sup\{a,b\}$ and $\inf\{a,b\}$, respectively, and writing $1_L$ and $0_L$, or simply $1$ and $0$ when no confusion will arise, respectively for the greatest element $\V L=\bigwedge\emptyset$ and the least element $\bigwedge L=\V\emptyset$ of a sup-lattice $L$.

As examples of sup-lattices we mention:
\begin{enumerate}
\item The topology $\topology(X)$ of a topological space $X$, ordered by inclusion; we have $\V S=\bigcup S$ for all $S\subset\topology(X)$, and $\bigwedge S$ is the interior of $\bigcap S$. For binary meets we have $U\wedge V=U\cap V$.
\item The set of linear subspaces, $\Sub A$, of a vector space $A$. If $(V_\alpha)$ is a family in $\Sub A$ we have $\V_\alpha V_\alpha =\spanmap \bigcup_\alpha V_\alpha$ and $\bigwedge_\alpha V_\alpha=\bigcap_\alpha V_\alpha$. Any linear map $f:A\to B$ between vector spaces yields a sup-lattice homomorphism $\Sub f:\Sub A\to\Sub B$ defined by
\[
\Sub f(V) = f(V)\;,
\]
and in this way one obtains a functor $\Sub : \Vect\to\SL$.
\item The set of closed linear subspaces, $\Max A$, of a topological vector space $A$. If $(V_\alpha)$ is a family in $\Max A$ we have $\V_\alpha V_\alpha =\overline{\spanmap\bigcup_\alpha V_\alpha}$ and $\bigwedge_\alpha V_\alpha=\bigcap_\alpha V_\alpha$. Any continuous linear map $f:A\to B$ between topological vector spaces yields a sup-lattice homomorphism $\Max f:\Max A\to\Max B$ defined by
\[
\Max f(V) = \overline{f(V)}\;,
\]
and in this way one obtains a functor $\Max : \TopVect\to\SL$.
\end{enumerate}

\paragraph{Adjunctions.}

Let $L$ and $M$ be sup-lattices. We write $f\dashv g$ in order to indicate that $f$ is left adjoint to $g$:
\[
\xymatrix{L\ar@/^2ex/[rr]^f&\perp&M\ar@/^2ex/[ll]^g}\;.
\]
The pair $(f,g)$ will be referred to as an \emph{adjunction from $L$ to $M$}.

\example
Let $A$ and $B$ be a complex vector spaces, and let $f:A\to B$ be a linear map. Recall that $\Sub f$ is the sup-lattice homomorphism $\Sub A\to\Sub B$ defined by $\Sub f(V)=f(V)$ for all $V\in\Sub A$. The right adjoint of $\Sub f$ is the inverse image map $f^{-1}:\Sub B\to\Sub A$, as the following equivalences show, where $V\in\Sub A$ and $W\in\Sub B$:
\[
\Sub f(V)\subset W\iff f(V)\subset W\iff V\subset f^{-1}(W)\;.
\]
\endexample

We note that the bijection between sup-lattice homomorphisms and their right adjoints is an antitone order isomorphism with respect to the pointwise order on maps: if both $f\dashv g$ and $f'\dashv g'$ are adjunctions from $L$ to $M$ then we have $f\le f'\iff g'\le g$.

\subsection{Locales}

\paragraph{Basic definitions and facts.}

A \emph{locale} (see \cite{stonespaces}) is a sup-lattice $L$ satisfying the distributivity property
\[
a\wedge\V_\alpha b_\alpha = \V_\alpha a\wedge b_\alpha
\]
for all $a\in L$ and all families $(b_\alpha)$ in $L$. The main example of a locale is the topology $\topology(X)$ of a topological space.

A \emph{map of locales} $f:M\to L$ is a homomorphism of sup-lattices in the opposite direction,
\[
f^*:L\to M\;,
\]
that also preserves finite meets:
\begin{eqnarray*}
f^*(1_L)&=& 1_M\\
f^*(a\wedge b) &=& f^*(a)\wedge f^*(b)\;.
\end{eqnarray*}
Then $f^*$ is called a \emph{homomorphism of locales}, and it is referred to as the \emph{inverse image homomorphism of the map $f$}. This terminology follows from the fact that every continuous map of topological spaces  $\varphi:Y\to X$ defines a map of locales $f:\topology(Y)\to\topology(X)$ by the condition:
\[
f^*=\varphi^{-1}:\topology(X)\to\topology(Y)\;.
\]

The \emph{category of locales} $\Loc$ has the locales as objects and the maps of locales as arrows. The assignment $X\mapsto\topology(X)$ extends to a functor
\[
\topology:\Top\to\Loc
\]
such that for each continuous map $\varphi$ we have
\[
(\topology \varphi)^*=\varphi^{-1}\;.
\]
We also note that, since
the inverse image homomorphism of a map of locales $f:M\to L$ preserves joins,  there is a meet preserving map (the right adjoint of $f^*$)
\[
f_*:M\to L\;,
\]
which is defined by, for all $b\in M$,
\[
f_*(b)=\V\{a\in L\st f^*(a)\le b\}\;.
\]
If $\varphi:Y\to X$ is a continuous open map between topological spaces then we denote by $\varphi_!:\topology(Y)\to\topology(X)$ the \emph{direct image homomorphism} defined by $\varphi_!(U)=\varphi(U)$. This is a sup-lattice homomorphism, left adjoint to $\varphi^{-1}:\topology(X)\to\topology(Y)$.

\paragraph{Points and prime elements.}

A \emph{point} of a locale $L$ is usually defined to be a map $x:\Omega\to L$, where $\Omega=\pwset{\{*\}}$ is the topology of a singleton. Although locale theory is often meant to be applied in a constructive setting (\ie, in an arbitrary topos), classically $\Omega$ is just a two element chain $\{0,1\}$ with $0<1$, and the set of points $x$ is in bijective correspondence with the set of elements of the form
\[
\V\ker x^* =\V\{a\in L\st x^*(a)=0\}\;.
\]
These elements are precisely the prime elements of $L$, where an element $p\in L$ of a locale $L$ is said to be \emph{prime} if it satisfies the following two conditions for all $a,b\in L$:
\begin{eqnarray*}
p&\neq& 1\\
a\wedge b\le p &\Rightarrow& a\le p\textrm{ or }b\le p\;.
\end{eqnarray*}
For the purposes of this paper it will be useful to identify the points of a locale $L$ with its prime elements, and we denote the set of the latter by $\spectrum(L)$.

The following simple fact will be used later:

\lemma\label{JPexl2}
Let $X$ be a topological space with a subbasis $\mathcal S$. Any prime element $P\in\Omega X$ is a union of subbasic open sets.
\endlemma
\proof
For any $x\in P$ there are subbasic open sets $S_1,\dots,S_k\in\mathcal S$
such that $x\in S_1\cap\dots\cap S_k\subset P$.
Since $P$ is prime, $x\in S_i\subset P$ for some $i$. 
\endproof

\paragraph{The spectrum of a locale.}

For each $a\in L$, where $L$ is a locale, we define the following set,
\[U_a=\{p\in \spectrum(L)\st a\nleq p\}\;,\]
and the following conditions hold:
\begin{eqnarray*}
U_1 &=& \spectrum(L)\;;\\
U_{a\wedge b} &=& U_a\cap U_b\;;\\
U_{\V_\alpha a_\alpha} &=&\bigcup_\alpha U_{a_\alpha}\;.
\end{eqnarray*}
Hence, the collection $(U_a)_{a\in L}$ defines a topology $\topology\spectrum(L)$ on $\spectrum(L)$. The set of primes equipped with this topology will be called the \emph{spectrum} of the locale $L$. The assignment $a\mapsto U_a$ is a homomorphism of locales that defines a map of locales $\spat_L:\topology\spectrum(L)\to L$ called the \emph{spatialization of $L$}. A locale is called \emph{spatial} if its spatialization map is an isomorphism. This is equivalent to the statement that $L$ is isomorphic to the topology of some topological space.

\paragraph{The spectrum functor.}

If $f:M\to L$ is a map of locales we obtain a continuous map
\[
\spectrum(f):\spectrum(M)\to\spectrum(L)
\]
which coincides with the restriction of $f_*$ to $\spectrum(M)$:
\[
\spectrum(f)(p) = \V\{a\in L\st f^*(a)\le p\}\;.
\]
In this way we have a functor
\[
\spectrum:\Loc\to\Top\;,
\]
which is right adjoint to $\topology$. If $X$ is a topological space there is a continuous map $\sob_X:X\to\spectrum\topology(X)$, called the \emph{soberification of $X$}, which to each $x\in X$ assigns the complement in $X$ of the closure $\overline{\{x\}}$. The space $X$ is $T_0$ if and only if $\sob_X$ is injective. Moreover, note the following simple fact:
\lemma\label{sobsurj}
If $\sob_X$ is surjective then it is an open map.
\endlemma
\proof
Let $W\subset X$ be open. The condition $\sob_X(x)\in U_W$ is equivalent to $x\in W$, so if $\sob_X$ is surjective we have $\sob_X(W) = U_W$. 
\endproof
It follows that if $\sob_X$ is bijective then it is necessarily a homeomorphism, and in this case $X$ is said to be a \emph{sober space}. In particular, the spectrum of a locale is a sober space, and hence $T_0$.

The adjunction between $\topology$ and $\spectrum$ restricts to an equivalence of categories between the full subcategories of spatial locales and sober spaces. If $X$ is a Hausdorff space, the prime elements of $\topology(X)$ are exactly the open sets of the form $X\setminus\{x\}$ and thus $X$ is sober.

\section{Quotient vector bundles via locales}\label{sec:qvbl}

The ideas described in the introduction show that a quotient vector bundle $\mathcal A=(\pi:E\to X,A,q)$ gives rise to an adjunction between $\Sub A$ and $\topology(X)$, both regarded as sup-lattices under the inclusion order. The left adjoint is the \emph{open support map} $\sigma:\Sub A\to\topology(X)$, and the right adjoint $\gamma:\topology(X)\to\Sub A$ yields, for each open set $U\subset X$, the set of ``formal sections'' $a\in A$ such that the open support $\osupp\hat a$ is contained in $U$. In this section we pursue this idea in order to study quotient vector bundles, in the end obtaining an adjunction between the category of such locales with linear structure and a subcategory of the category of quotient vector bundles. The objects of this category will be called \emph{spectral vector bundles}. Any quotient vector bundle on a sober base space and with closed zero section will be seen to be of this kind.

\subsection{Linearized locales}

\paragraph{Basic definitions and facts.}

Let $\diag$ be a locale. By a \emph{linear structure} on $\diag$ is meant a topological vector space $A$ together with a sup-lattice homomorphism
\[
\sigma:\Sub A\to\diag\;,
\]
which we refer to as the \emph{support map}, and whose right adjoint $\gamma:\diag\to\Sub A$ restricts to a continuous map $\kfrak=\gamma\vert_{\spectrum(\diag)}:\spectrum(\diag)\to\Sub A$. The map $\gamma$ is referred to as the \emph{restriction map} --- for each $U\in\diag$ we may think of $\gamma(U)$ as the restriction of $A$ to $U$, where $A$ is regarded as a space of formal global sections on $\diag$. The map $\kfrak$ is referred to as the \emph{kernel map} of the linear structure.

By a \emph{linearized locale} $\mathfrak A=(\diag,A,\sigma,\gamma)$ will be meant a locale $\diag$ equipped with a linear structure given by $A$ and $\sigma$, with $\sigma\dashv \gamma$.

The continuity of the kernel map of a linear structure has the following necessary condition:

\lemma\label{kernelclassfunction}
Let $\mathfrak A=(\diag,A,\sigma,\gamma)$ be a linearized locale. For all $p,q\in\spectrum(\diag)$ the following implication holds:
\[p\le q\Rightarrow \overline{\gamma(p)}=\overline{\gamma(q)}\;.\]
\endlemma

\proof
The specialization order of the topology of $\spectrum(\diag)$ is dual to the order of $\diag$. Hence, the kernel map, which is monotone because $\gamma$ is, is also antitone on prime elements due to continuity, and thus the condition $p\le q$ implies that $\gamma(p)$ and $\gamma(q)$ are topologically equivalent, which for the lower Vietoris topology means precisely that their closures in $A$ are the same. 
\endproof

We also mention, although we shall not use it in this paper, the following general necessary and sufficient condition for the continuity of the kernel map.

\theorem
Let $\diag$ be a locale, $A$ a topological vector space, and $\sigma\dashv\gamma$ an adjunction from $\Sub A$ to $\diag$. Then $\kfrak:=\gamma\vert_{\spectrum(\diag)}$ is continuous if and only if
for all open sets $U\subset A$ the set $\bigcap_{a\in U} U_{\sigma(\linspan a)}$ is closed in $\spectrum(\diag)$.
\endtheorem

\proof
Let $U\subset A$ be an open set and $p\in\spectrum(\diag)$. Then
\begin{eqnarray*}
p\in \kfrak^{-1}(\tilde U) &\iff& \kfrak(p)\cap U\neq\emptyset\\
&\iff& \exists_{a\in U}\ a\in\kfrak(p)\\
&\iff& \exists_{a\in U}\ \linspan a\subset\gamma(p)\\
&\iff& \exists_{a\in U}\ \sigma(\linspan a)\le p\\
&\iff& \exists_{a\in U}\ p\notin U_{\sigma(\linspan a)}\\
&\iff& p\notin\bigcap_{a\in U} U_{\sigma(\linspan a)}\;,
\end{eqnarray*}
and thus
$
\kfrak^{-1}(\tilde U) = \spectrum(\diag)\setminus\bigcap_{a\in U} U_{\sigma(\linspan a)}
$.
Hence, $\kfrak$ is continuous if and only if $\bigcap_{a\in U} U_{\sigma(\linspan a)}$ is closed for all $U\in\topology(A)$. 
\endproof

\paragraph{Morphisms of linearized locales.}

Let the following be linearized locales:
\begin{eqnarray*}
\mathfrak A&=&(A,\diag_A,\sigma_A,\gamma_A)\;,\\
\mathfrak B&=&(B,\diag_B,\sigma_B,\gamma_B)\;.
\end{eqnarray*}
A \emph{morphism} $\mathfrak f:\mathfrak B\to\mathfrak A$ is a pair $(\underline f,\overline f)$ consisting of a map of locales
\[\underline f:\diag_B\to\diag_A\]
and a continuous linear map
\[\overline f:A\to B\]
satisfying, for all $V\in\Sub A$, the inclusion
\begin{equation}\label{llmorphism1}
\sigma_B(\overline f(V))\subset\underline f^*(\sigma_A(V))\;.
\end{equation}
In other words, $\mathfrak f$ satisfies the following lax commutation relation:
\begin{equation}\label{llmorphism2}
\vcenter{\xymatrix{
\Sub A\ar[d]_{\sigma_A}\ar[rr]^{\Sub \overline f}&~\ar@{}[d]|{\ge}&\Sub B\ar[d]^{\sigma_B}\\
\diag_A\ar[rr]_{\underline f^*}&~&\diag_B
}}
\end{equation}
Equivalently, on right adjoints, recalling that the right adjoint of $\Sub \overline f$ coincides with $\overline f^{-1}$, we have:
\begin{equation}\label{llmorphism3}
\vcenter{\xymatrix{
\Sub A&~\ar@{}[d]|{\le}&\Sub B\ar[ll]_{\overline f^{-1}}\\
\diag_A\ar[u]^{\gamma_A}&~&\diag_B\ar[u]_{\gamma_B}\ar[ll]^{\underline f_*}
}}
\end{equation}

These morphisms yield an obvious category $\LLoc$, which we shall refer to as the \emph{category of  linearized locales (with lax morphisms)}. By requiring the above commutation relations to be strict we obtain the subcategory $\sLLoc$, referred to as the \emph{strict category of linearized locales}. Natural examples of such morphisms can be obtained from Example~\ref{examplemorphQVBun} via the adjunction that will be obtained later in this section.

\subsection{Spectral vector bundles}

\paragraph{The spectral kernel of a quotient vector bundle.}

Let $\mathcal A=(\pi:E\to X,A,q)$ be a quotient vector bundle with kernel map $\kappa:X\to\Sub A$, and consider the adjoint pair
\[
\xymatrix{
\Sub A\ar@/^2ex/[rr]^-{\sigma}&\perp&\diag\ar@/^2ex/[ll]^-{\gamma}
}
\]
where $\sigma$ and $\gamma$ are defined by
\begin{eqnarray*}
\sigma(V)&=&\bigcup_{a\in V}\osupp\hat a\;,\\
\gamma(U)&=&\spanmap\{a\in A\st \osupp \hat a\subset U\}\;.
\end{eqnarray*}
The restriction of $\gamma$ to the set of prime open sets of $X$,
\[
\kfrak:=\gamma\vert_{\spectrum(\diag)}:\spectrum(\diag)\to\Sub A\;,
\]
will be referred to as the \emph{spectral kernel map} of the bundle $\mathcal A$. Then we have the following immediate remark:

\lemma\label{linlociffkfrakcont}
$(\topology(X),A,\sigma,\gamma)$ is a linearized locale if and only if the spectral kernel $\kfrak$ of $\mathcal A$ is continuous.
\endlemma

Also, note the following relation between the two kernel maps:

\lemma\label{laxvsstrict0}
For all $x\in X$ we have $\kappa(x)\subset\kfrak\bigl(\sob x\bigr)$.
\endlemma

\proof
We have $\kappa(x)=\{a\in A\st \hat a(x)=0\}$ and
$\kfrak(\sob x)=\kfrak(X\setminus\overline{\{x\}})=\{a\in A\st x\notin\osupp\hat a\}$,
and also $\hat a(x)=0\Rightarrow x\notin\osupp\hat a$. 
\endproof

\paragraph{Open support property.}

In most examples that arise in practice, such as those obtained from Banach bundles, the (image of the) zero section of a quotient vector bundle is closed in $E$, and thus the sets $\{x\in X\st \hat a(x)\neq 0\}$ are open, so that we have
\[
\osupp \hat a = \{x\in X\st \hat a(x)\neq 0\}\;.
\]
We will say that the quotient vector bundle has the \emph{open support property} if the above condition holds for all $a\in A$. As will be seen below, this property is in general weaker than having a closed zero section and it is useful in its own right.

The open support property can also be formulated equivalently in terms of the kernel maps $\kappa$ and $\kfrak$ of $\mathcal A$, as we now show.

\lemma\label{laxvsstrict}
$\mathcal A$ has the open support property if and only if $\kappa = \kfrak\circ\sob$. In particular, these equivalent conditions hold whenever the (image of the) zero section of $\pi$ is a closed set of $E$.
\endlemma

\proof
We have $\kappa=\kfrak\circ\sob$ if and only if for all $a\in A$ and all $x\in X$ the following equivalence holds:
\[\hat a(x)=0\iff x\notin\osupp\hat a\;. \] 
\endproof

\paragraph{Continuous spectral kernels.}

We say that $\mathcal A$ is a \emph{spectral (quotient) vector bundle} if $\mathcal A$ has the open support property and the spectral kernel $\kfrak$ is continuous. 
The full subcategory of $\QVBun$ (resp.\ $\sQVBun$) whose objects are the spectral vector bundles is denoted by $\pQVBun$ (resp.\ $\spQVBun$).

The following is a useful sufficient condition for spectrality:

\lemma\label{surjsob}
If a quotient vector bundle $\mathcal A=(\pi:E\to X,A,q)$ has the open support property and $\sob:X\to\spectrum\topology(X)$ is surjective then $\mathcal A$ is spectral.
\endlemma

\proof
The surjectivity of $\sob$ tells us that the direct image $\sob_!:\topology(X)\to\topology\spectrum\topology(X)$ is surjective (\cf\ Lemma~\ref{sobsurj}). Therefore the inverse image map $\sob^{-1}$ is injective and
$\sob_!\circ\sob^{-1}$ is the identity map on $\topology\spectrum\topology(X)$.
Hence,
\[
\kfrak^{-1}=\sob_!\circ\sob^{-1}\circ\kfrak^{-1}=\sob_!\circ\kappa^{-1}
\]
and it follows that $\kfrak$ is continuous. 
\endproof

\example\label{banachspectralexm}
Any Banach bundle $\pi:E\to X$ as in \cite{FD1} has a Hausdorff base space $X$ and a closed zero section, and thus any quotient vector bundle $(\pi,A,q)$ is spectral. For instance, if $X$ is locally compact we may take $A=C_0(\pi)$ and $q=\eval$ \cite{RS}.
\endexample

\example
It is possible to have a continuous spectral kernel and a sober base space without having the open support property. An example of this is the universal quotient vector bundle with Hausdorff fibers for a locally convex space $A$, whose base space $\Max A$ is sober --- see Theorem~\ref{JPext1} below.
\endexample

\example
For examples with the open support property but discontinuous spectral kernel see Theorem~\ref{JPext2} below, where it is also shown that for suitable $A$ (in particular finite dimensional) we obtain spectral vector bundles whose base spaces are sober and $T_1$ (\cf\ Lemma~\ref{oMaxT1}), but not Hausdorff (\cf\ Example~\ref{oMaxnotHaus}).
\endexample

There is also a necessary condition for spectrality in terms of the kernel map $\kappa$, based on Lemma~\ref{kernelclassfunction}:

\lemma\label{kernelclassfunctionbundles}
Let $\mathcal A = (\pi:E\to X,A,q)$ be a spectral vector bundle with kernel map $\kappa$. For all $x,y\in X$ we have the following implication, where $\sqsubseteq$ denotes the specialization preorder of the topology of $X$:
\[
x\sqsubseteq y\Rightarrow\overline{\kappa(x)}=\overline{\kappa(y)}\;.
\]
\endlemma

\proof
Let $x,y\in X$ be such that $x\sqsubseteq y$. Then $\sob(x)\sqsubseteq\sob(y)$ in $\spectrum(\diag)$, which is equivalent to $\sob(y)\le\sob(x)$, and thus by Lemma~\ref{kernelclassfunction} we obtain
\[
\overline{\kappa(x)} = \overline{\gamma(\sob(x))}=\overline{\gamma(\sob(y))}=\overline{\kappa(y)}\;. 
\]
\endproof

Less obvious is that Lemma~\ref{kernelclassfunction} has a converse for bundles, which gives us a necessary and sufficient condition for spectrality:

\theorem\label{thmcontkfrak}
Let $\mathcal A=(\pi:E\to X,A,q)$ be a quotient vector bundle satisfying the open support property. Then the spectral kernel $\kfrak$ is continuous (and thus $\mathcal A$ is a spectral vector bundle) if and only if for all $P,Q\in \spectrum\topology(X)$ we have the implication
\begin{equation}\label{thmimplication}
P\subset Q\Rightarrow\overline{\kfrak(P)}=\overline{\kfrak(Q)}\;.
\end{equation}
\endtheorem

\proof
If $\kfrak$ is continuous we obtain \eqref{thmimplication} from Lemma~\ref{kernelclassfunction} and Lemma~\ref{linlociffkfrakcont}, so
it remains to be proved only that \eqref{thmimplication} implies the continuity of $\kfrak$. In order to prove this we shall show that for any open set $W\subset\Sub A$ we have
\begin{equation}\label{auxeq}
\kfrak^{-1}(W)=U_{\kappa^{-1}(W)}\;,
\end{equation}
where $\kappa$ is the kernel map.
Let $P\in U_{\kappa^{-1}(W)}$. Then $\kappa^{-1}(W)\not\subset P$, so there is
$x\in X$ such that $x\in \kappa^{-1}(W)$ and $x\notin P$. But
then $\kappa(x)\in W$ and
\[
x\notin P\ \iff\ P\subset\sob(x)\ \Longrightarrow\ 
\kfrak(P)\equiv\kfrak\bigl(\sob(x)\bigr)=\kappa(x)\;,
\]
where $\equiv$ denotes topological equivalence.
Since $\kappa(x)\in W$ and $\kfrak(P)\equiv \kappa(x)$ we have $\kfrak(P)\in W$, and thus
$P\in\kfrak^{-1}(W)$.

Reciprocally, assume $P\in\kfrak^{-1}(W)$, that is, 
$\kfrak(P)\in W$. Then for any $x\notin P$
we have $P\subset \sob(x)$, so 
\[
\kappa(x)=\kfrak\bigl(\sob(x)\bigr)\equiv\kfrak(P)\;,
\]
and hence $\kappa(x)\in W$.
This shows that $X\setminus P\subset \kappa^{-1}(W)$. But then,
since $P\neq X$, we have $\kappa^{-1}(W)\not\subset P$, so
$P\in U_{\kappa^{-1}(W)}$, and this proves \eqref{auxeq}. 
\endproof

\paragraph{The spectrum of a linearized locale.}

Let $\mathfrak A=(\diag,A,\sigma,\gamma)$ be a linearized locale. We define its \emph{spectrum} to be the quotient vector bundle
\[
\specfunctor(\mathfrak A) = \kfrak^*(UA)
\]
which is classified by the continuous map
\[
\kfrak=\gamma\vert_{\spectrum(\diag)}:\spectrum(\diag)\to\Sub A\;.
\]

\lemma\label{specisspectral}
For any linearized locale $\mathfrak A$, the spectrum $\specfunctor(\mathfrak A)$ is a spectral vector bundle.
\endlemma

\proof
Let us use the following notation:
\[
\bunfunctor\specfunctor(\mathfrak A) = (\topology\spectrum(\diag),A,\sigma,\gamma)\;.
\]
The base space of $\specfunctor(\mathfrak A)$ is $\spectrum(\diag)$, which is a sober space, so, by Lemma~\ref{laxvsstrict} and Lemma~\ref{surjsob}, we need only verify that for each $a\in A$ the set
\[
\{p\in\spectrum(\diag)\st a\notin\gamma(p)\}
\]
is open in $\spectrum(\diag)$. This follows from the adjunction $\sigma\dashv\gamma$, because for all $p\in\spectrum(\diag)$ and $a\in A$ we have
\[
\sigma(\linspan a)\le p\iff \linspan a\subset\gamma(p)\;,
\]
and thus $\{p\in\spectrum(\diag)\st a\notin\gamma(p)\}$ is the open set $U_{\sigma(\linspan a)}$. 
\endproof

\theorem
The assignment $\mathfrak A\mapsto\specfunctor(\mathfrak A)$ extends to a functor
\[
\specfunctor:\LLoc\to\pQVBun\;.
\]
Moreover, this functor restricts to the corresponding strict categories:
\[
\specfunctor:\sLLoc\to\spQVBun\;.
\]
\endtheorem

\proof
Let $\mathfrak A=(\diag,A,\sigma,\gamma)$ and $\mathfrak A'=(\diag',A',\sigma',\gamma')$ be linearized locales, and let
\[
\mathfrak f=(\underline f,\overline f):\mathfrak A'\to\mathfrak A
\]
be a morphism in $\LLoc$. Defining
\begin{eqnarray*}
f_\flat&=&\underline f_*\vert_{\spectrum(\diag')}:\spectrum(\diag')\to\spectrum(\diag)\\
f^*&=&\overline f:A\to A'
\end{eqnarray*}
we obtain a morphism $f:\specfunctor(\mathfrak A')\to\specfunctor(\mathfrak A)$ in $\QVBun$. In order to see this we only need to prove that \eqref{qvbmorphism} holds. Let $a\in A$ and $y\in\spectrum(\diag')$, and denote by $\kfrak$ and $\kfrak'$ the restrictions of $\gamma$ and $\gamma'$, respectively, to $\spectrum(\diag)$ and $\spectrum(\diag')$. Denote also the quotient vector bundles as follows:
\begin{eqnarray*}
\specfunctor(\mathfrak A) &=& (\pi:E\to\spectrum(\diag),A,q)\\
\specfunctor(\mathfrak A') &=& (\pi':E'\to\spectrum(\diag'),A',q')
\end{eqnarray*}
Then, for all $a\in A$ and $y\in\spectrum(\diag')$, we have
\begin{eqnarray*}
q(a,f_\flat(y))=0 &\iff& a\in\kfrak(f_\flat(y))\\
&\iff&a\in\gamma\circ \underline f_*(y)\\
&\Longrightarrow& a\in\overline f^{-1}\circ\gamma'(y)\\
&\iff& \overline f(a)\in\kfrak'(y)\\
&\iff& q'(f^*(a),y)=0\;.
\end{eqnarray*}
Hence, \eqref{qvbmorphism} holds. It is also clear that if $\mathfrak f$ is a strict morphism then the above implication is another equivalence, and thus $f$ is also strict. The functoriality is evident, so we have obtained two functors as intended. 
\endproof

\subsection{Adjunctions and equivalences}

\paragraph{A left adjoint to the spectrum functor.}

For each spectral vector bundle
\[
\mathcal A=(\pi:E\to X,A,q)\;,
\]
let
\[\bunfunctor(\mathcal A)=(\topology(X),A,\sigma,\gamma)\] be the corresponding linearized locale.

\lemma
Let $\mathcal A=(\pi:E\to X,A,q)$ be a spectral vector bundle with kernel map $\kappa$, and let
\[
\specfunctor\bunfunctor(\mathcal A)=(\mathfrak p:\mathfrak E\to\spectrum\topology(X),A,\mathfrak q)\;,
\]
with kernel map $\kfrak$.
Then
\[\sob_{\mathcal A}:=(\sob_X,\ident_A):\mathcal A\to\specfunctor\bunfunctor(\mathcal A)\]
is a morphism in $\spQVBun$.
\endlemma

\proof
The condition for being a strict morphism is given by \eqref{qvbmorphism3}, so for all $a\in A$ and $x\in X$ we must prove the equivalence
\[
\mathfrak q(a,\sob_X(x))=0\iff q(a,x)=0\;,
\]
which is immediate because $\kfrak(\sob_X(x))=\kappa(x)$:
\[\mathfrak q(a,\sob_X(x))=0\iff a\in\kfrak(\sob_X(x))\iff a\in\kappa(x)\iff q(a,x)=0\;.\]
Hence, since both $\mathcal A$ and  $\specfunctor\bunfunctor(\mathcal A)$ are spectral  (\cf\ Lemma~\ref{specisspectral}), $\sob_{\mathcal A}$ is indeed a morphism in
$\spQVBun$. 
\endproof

\remark
Note that for $\sob_{\mathcal A}$ to be a morphism at all it is necessary that it be strict, for the weaker condition \eqref{qvbmorphism} in this case reads
\[\mathfrak q(a,\sob_X(x))=0\quad\Longrightarrow\quad q(a,x)=0\;,\]
which implies
\[\kfrak\circ\sob_X\le\kappa\]
and, by Lemma~\ref{laxvsstrict}, is equivalent to strictness. Hence, for any quotient vector bundle $\mathcal A$ with continuous $\kfrak$, the following conditions are equivalent:
\begin{enumerate}
\item $\mathcal A$ is spectral;
\item $\sob_{\mathcal A}$ is a morphism;
\item $\sob_{\mathcal A}$ is a strict morphism.
\end{enumerate}
\endremark

\theorem
The assignment $\mathcal A\mapsto\bunfunctor(\mathcal A)$ extends to a functor
\[\bunfunctor:\pQVBun\to\LLoc\;,\]
which is left adjoint to $\specfunctor$. Moreover the unit of the adjunction is the family $\{\sob_{\mathcal A}\}_{\mathcal A}$, and the adjunction restricts to the strict subcategories $\spQVBun$ and $\sLLoc$.
\endtheorem

\proof
In order to prove the theorem we shall show, for an arbitrary spectral vector bundle $\mathcal A$, that $(\bunfunctor(\mathcal A),\sob_{\mathcal A})$ is a universal arrow from $\mathcal A$ to $\specfunctor$.
So let the following be, respectively, a spectral vector bundle and a linearized locale,
\begin{eqnarray*}
\mathcal A&=&(\pi:E\to X,A,q)\;,\\
\mathfrak B&=&(\diag,B,\sigma_B,\gamma_B)\;,
\end{eqnarray*}
and let $f=(f_\flat,f^*):\mathcal A\to\specfunctor(\mathfrak B)$ be a morphism in $\pQVBun$ (\cf\ Lemma~\ref{specisspectral}). We need to prove that there is a unique morphism of linearized locales
\[
\mathfrak f=(\underline f,\overline f):\bunfunctor(\mathcal A)\to\mathfrak B
\]
such that the following diagram commutes:
\[
\xymatrix{
\mathcal A\ar[rr]^{\sob_{\mathcal A}}\ar[drr]_{f}&&\specfunctor\bunfunctor(\mathcal A)\ar[d]^{\specfunctor(\mathfrak f)}\\
&&\specfunctor(\mathfrak B)
}
\]
Let us use the following notations for the structure maps of $\bunfunctor(\mathcal A)$ and $\specfunctor(\mathfrak B)$:
\begin{eqnarray*}
\bunfunctor(\mathcal A)&=&(\topology(X),A,\sigma_A,\gamma_A)\;,\\
\specfunctor(\mathfrak B)&=&(\pi': E'\to\spectrum(\diag),B, q')\;.
\end{eqnarray*}
Due to the adjunction $\topology\vdash\spectrum$ between $\Top$ and $\Loc$ there is in $\Loc$ a unique map $\underline f:\topology(X)\to\diag$ such that $\spectrum(\underline f)\circ\sob_X=f_\flat$ --- its inverse image is defined, for all $d\in\diag$, by
\begin{equation}\label{llmorphismadj1}
\underline f^*(d) = f_\flat^{-1}(U_d)\;,
\end{equation}
where $U_d=\{p\in\spectrum(\diag)\st d\nleq p\}$.
In addition, $\overline f$ must obviously be $f^*$, so we must now verify that \eqref{llmorphism1} holds for $\mathfrak f$. With the necessary notation changes this condition is, for all $V\in\Sub B$,
\begin{equation}\label{llmorphismadj2}
\sigma_A(f^*(V))\subset\underline f^*(\sigma_B(V)) \;,
\end{equation}
or, using \eqref{llmorphismadj1},
\begin{equation}\label{llmorphismadj3}
\sigma_A(f^*(V))\subset f_\flat^{-1}(U_{\sigma_B(V)})\;.
\end{equation}
In order to prove the latter condition we shall prove that, for all $x\in X$, if $x\notin f_\flat^{-1}(U_{\sigma_B(V)})$ then $x\notin \sigma_A(f^*(V))$. Let then $x\notin f_\flat^{-1}(U_{\sigma_B(V)})$. This is equivalent to $f_\flat(x)\notin U_{\sigma_B(V)}$, which in turn is equivalent to $\sigma_B(V)\leq f_\flat(x)$, and, by adjointness, to
\begin{equation}\label{llmorphismadj4}
V\subset\gamma_B(f_\flat(x))\;.
\end{equation}
Now $f_\flat(x)\in\spectrum(\diag)$ and thus, by definition, $\gamma_B(f_\flat(x))= \kappa'(f_\flat(x))$, where $\kappa'$ is the kernel map of $\specfunctor(\mathfrak B)$. Hence, \eqref{llmorphismadj4} is equivalent to the statement that every $a\in V$ satisfies $a\in \kappa'(f_\flat(x))$, which in turn is equivalent to
\begin{equation}\label{llmorphismadj5}
\forall_{a\in V}\ q'(a,f_\flat(x))= 0\;.
\end{equation}
Since $f$ is a morphism of quotient vector bundles the latter implies
\begin{equation}\label{llmorphismadj6}
\forall_{a\in V}\ q(f^*(a),x)= 0\;,
\end{equation}
which is equivalent to
\begin{equation}\label{llmorphismadj7}
f^*(V)\subset\kappa(x)\;.
\end{equation}
Since $\mathcal A$ is spectral we have $\kappa=\kfrak\circ\sob_X$, and thus \eqref{llmorphismadj7} is equivalent to
\[
f^*(V)\subset\kfrak(\sob_X(x))=\gamma_A(\sob_X(x))\;,
\]
and, by adjointness, equivalent to
\[
\sigma_A(f^*(V))\subset\sob_X(x)=X\setminus\overline{\{x\}}\;,
\]
which, since $\sigma_A(f^*(V))$ is open in $X$, is equivalent to
\begin{equation}\label{llmorphismadj8}
x\notin \sigma_A(f^*(V))\;.
\end{equation}
This concludes the proof of the inclusion \eqref{llmorphismadj3}, and establishes the adjunction $\bunfunctor\vdash\specfunctor$. 
Now we note that all but one of the steps in the above proof of \eqref{llmorphismadj2} are equivalences. But
if $f$ is a strict morphism the implication \eqref{llmorphismadj5}$\Rightarrow$\eqref{llmorphismadj6} becomes an equivalence, too, and
we obtain
\[
\sigma_A(f^*(V))=\underline f^*(\sigma_B(V))\;,
\]
which means that $\mathfrak f$ is strict. Together with the fact that $\sob_{\mathcal A}$ is strict, this proves that the adjunction restricts to the strict subcategories. 
\endproof

\paragraph{Sober vector bundles.} Let $\mathcal A=(\pi:E\to X,A,q)$ be a quotient vector bundle with kernel map $\kappa:X\to\Sub A$. We say that $\mathcal A$ is a \emph{sober vector bundle} if the two following conditions hold:
\begin{itemize}
\item $X$ is a sober space;
\item $\mathcal A$ has the open support property.
\end{itemize}

By Lemma~\ref{laxvsstrict} and Lemma~\ref{surjsob} we immediately obtain:

\begin{corollary}
The following conditions satisfy $\eqref{cond1}\Rightarrow \eqref{cond2}\Rightarrow \eqref{cond3}$:
\begin{enumerate}
\item\label{cond1} $X$ is sober and the zero-section of $\mathcal A$ is closed in $E$;
\item\label{cond2} $\mathcal A$ is sober;
\item\label{cond3} $\mathcal A$ is spectral.
\end{enumerate}
\end{corollary}

\example
Any Banach bundle on a locally compact Hausdorff space can be made a sober vector bundle (\cf\ Example~\ref{banachspectralexm}). Examples of quotient vector bundles that are spectral but not sober can of course be obtained from any sober vector bundle $(\pi:E\to X,A,q)$ by pulling it back along a quotient map $Y\to X$ such that $Y$ is not sober.
\endexample

\theorem
Let $\mathfrak A=(\diag,A,\sigma,\gamma)$ and $\mathcal A=(\pi,A,q)$ be a linearized locale and a quotient vector bundle, respectively.
\begin{enumerate}
\item\label{AA1} The spectrum $\specfunctor(\mathfrak A)$ is a sober vector bundle.
\item\label{AA2} The bundle $\mathcal A$ is sober if and only if
\[
\sob_{\mathcal A}=(\sob_X,\ident_A):\mathcal A\to\specfunctor\bunfunctor(\mathcal A)
\]
is an isomorphism.
\end{enumerate}
\endtheorem

\proof
We have already seen in Lemma~\ref{specisspectral} that $\specfunctor(\mathfrak A)$ is spectral, so it is sober because $\spectrum(\diag)$ is. This proves \eqref{AA1}, and \eqref{AA2} is obvious. 
\endproof

\paragraph{Spatial linearized locales.} Let $\mathfrak A=(\diag,A,\sigma,\gamma)$ be a linearized locale. The \emph{spatialization} of $\mathfrak A$ is the pair
\[
\spat_{\mathfrak A} = (\spat_{\diag}:\topology\spectrum(\diag)\to\diag,\ident_A:A\to A)\;.
\]
\theorem
$\spat_{\mathfrak A}$ is a strict morphism of linearized locales
$\bunfunctor\specfunctor(\mathfrak A)\to\mathfrak A$.
\endtheorem

\proof
Let $\bunfunctor\specfunctor(\mathfrak A)=(\topology\spectrum(\diag),A,\tilde\sigma,\tilde\gamma)$. Then $\spat_{\mathfrak A}$ is a strict morphism if and only if the following square commutes:
\[
\xymatrix{
\Sub A\ar[rr]^{\ident_A}\ar[d]_{\sigma}&&\Sub A\ar[d]^{\tilde\sigma}\\
\diag\ar[rr]_{\spat_{\diag}^*}&&\topology\spectrum(\diag)
}
\]
Let $V\in\Sub A$. We have
\begin{eqnarray*}
\tilde\sigma(V) &=& \bigcup_{a\in V} \interior\{p\in\spectrum(\diag)\st a\notin\gamma(p)\}\\
&=& \bigcup_{a\in V} \{p\in\spectrum(\diag)\st a\notin\gamma(p)\}\quad\textrm{(by Lemma~\ref{laxvsstrict} and Lemma~\ref{specisspectral})}\\
&=&\bigcup_{a\in V} \{p\in\spectrum(\diag)\st \sigma(\linspan a)\nleq p\}\quad\textrm{(because $\sigma\dashv\gamma$)}\\
&=&\bigcup_{a\in V} U_{\sigma(\linspan a)}= U_{\sigma(V)}=\spat_{\diag}^*(\sigma(V))\;,
\end{eqnarray*}
and thus we have $\tilde\sigma = {\spat_{\diag}^*}\circ{\sigma}$, as intended. 
\endproof

We say that $\mathfrak A$ is \emph{spatial} if $\diag$ is a spatial locale. This is equivalent to stating that $\spat_{\mathfrak A}$ is an isomorphism of linearized locales. It immediately follows that the adjunction $\bunfunctor\dashv\specfunctor$ restricts to another adjunction, between the full subcategories of sober vector bundles and spatial linearized locales, whose unit and co-unit are isomorphisms. Hence, we are provided with a linearized version of the equivalence of categories between sober spaces and spatial locales:

\begin{corollary}
The adjunction $\bunfunctor\dashv\specfunctor$ restricts to an equivalence between the full subcategories of $\pQVBun$ and $\LLoc$ whose objects are, respectively, the sober vector bundles and the spatial linearized locales. This further restricts to an equivalence between the respective subcategories whose morphisms are strict.
\end{corollary}

\subsection{Bundles with Hausdorff fibers}

\paragraph{$\Max$-valued linear structures.} A quotient vector bundle $(\pi,A,q)$ has Hausdorff fibers if and only if its kernel map is valued in $\Max A$. All our previous results carry through to this setting, and we briefly examine this.

First, we notice that $\Max A$ is both a topological retract of $\Sub A$ (under the lower Vietoris topology) and a quotient sup-lattice given by the topological closure operator. In other words, we have an adjunction
\begin{equation}\label{submaxadj}
\xymatrix{\Sub A\ar@/^2ex/[rr]^{\overline{(-)}}&\perp&\Max A\ar@/^2ex/[ll]^\iota}
\end{equation}
where both the inclusion $\iota$ and the surjection $\overline{(-)}$ are continuous.

Now let $\mathfrak A=(\diag,A,\sigma,\gamma)$ be a linearized locale. Then $\gamma$ is valued in $\Max A$ if and only if $\sigma$ factors through $\Max A$ --- that is, $\sigma(V)=\sigma(\overline V)$ for all $V\in\Sub A$ ---, in which case we obtain an adjunction
\[
\xymatrix{\Max A\ar@/^2ex/[rr]^-{\sigma\vert_{\Max A}}&\perp&\diag\ar@/^2ex/[ll]^-\gamma}\;,
\]
and the kernel map $\kfrak$ is continuous as a map to $\Max A$.

Conversely, given an adjunction
\begin{equation}\label{maxlinloc}
\xymatrix{\Max A\ar@/^2ex/[rr]^-{\sigma}&\perp&\diag\ar@/^2ex/[ll]^-\gamma}
\end{equation}
such that $\kfrak:=\gamma\vert_{\spectrum(\diag)}:\spectrum(\diag)\to\Max A$ is continuous, composing this with the adjunction of \eqref{submaxadj} we obtain a linearized locale $(\diag,A,\overline{(-)}\circ\sigma,\iota\circ\gamma)$, whose kernel map is $\iota\circ\kfrak$.

Let us refer to such a linearized locale as a \emph{$\Max$-linearized locale}. All our results in section \ref{sec:qvbl} remain true if restricted to $\Max$-linearized locales and spectral vector bundles with Hausdorff fibers. In addition
we note that, by Lemma~\ref{kernelclassfunction}, if $(\diag,A,\sigma,\gamma)$ is a $\Max$-linearized locale then for all $p,q\in\spectrum(\diag)$ we have the following implication:
\begin{equation}\label{Maxvalkernelclassfunction}
p\le q\Rightarrow\gamma(p)=\gamma(q)\;.
\end{equation}

\paragraph{Spectral vector bundles revisited.}
Lemma~\ref{kernelclassfunctionbundles} has the following simple consequences:

\begin{corollary}
Let $\mathcal A = (\pi:E\to X,A,q)$ be a spectral vector bundle with Hausdorff fibers and kernel map $\kappa$. For all $x,y\in X$ we have the following implication, where $\sqsubseteq$ denotes the specialization preorder of the topology of $X$:
\[
x\sqsubseteq y\Rightarrow\kappa(x)=\kappa(y)\;.
\]
\end{corollary}

\begin{corollary}
Spectral vector bundles with Hausdorff fibers and injective kernel maps must have $T_1$ base spaces.
\end{corollary}

We also rewrite Theorem~\ref{thmcontkfrak} for bundles with Hausdorff filbers:

\begin{corollary}
Let $\mathcal A=(\pi:E\to X,A,q)$ be a quotient vector bundle with Hausdorff fibers satisfying the open support property.
Then the spectral kernel $\kfrak$ is continuous (and thus $\mathcal A$ is a spectral vector bundle) if and only if for all $P,Q\in \spectrum\topology(X)$ we have the implication
\begin{equation}\label{thmimplication2}
P\subset Q\Rightarrow\kfrak(P)=\kfrak(Q)\;.
\end{equation}
\end{corollary}

\section{Classifying spaces and universal bundles}

In this section we define the \emph{open support topology} on the spectrum $\Sub A$ of an arbitrary topological vector space $A$. This topology is finer than the lower Vietoris topology (for instance $\Max A$ becomes a $T_1$ space) but in general coarser than the Fell topology. We write $\oSub A$ for $\Sub A$ with the open support topology. This space classifies the spectral vector bundles, at least those over base spaces $X$ such that $\sob_X$ is surjective (\eg, Hausdorff spaces). We also show that this classification does not necessarily apply to more general spaces $X$ because although the quotient vector bundles classified by continuous maps $\kappa:X\to\Sub A$ have the open support property, their spectral kernels may fail to be continuous. In particular, in general there is no universal spectral vector bundle for a given topological vector space $A$. Conversely, we prove, for any locally convex space $A$, that the quotient vector bundle determined by the inclusion $\Max A\to\Sub A$ (with the lower Vietoris topology) does not have the open support property but that its spectral kernel is continuous (in fact constant). Along the way we obtain general results concerning the prime open sets of topologies on $\Sub A$ that are finer than the lower Vietoris topology, in particular concluding that for locally convex spaces $A$ the space $\Max A$ is sober.

\subsection{Classifying spaces for spectral vector bundles}

\paragraph{Open support topology.}
Let $A$ be a topological vector space. For each $a\in A$ let
\[
\check a = \{P\in\Sub A\st a\notin P\}\;.
\]
The coarsest topology on $\Sub A$ that contains the lower Vietoris topology and makes all the sets $\check a$ open will be referred to as the \emph{open support topology}.
The terminology is motivated by the following lemma:

\lemma
Let $\mathcal A$ be a quotient vector bundle. Then $\mathcal A$ has the open support property if and only if its kernel map is continuous with respect to the open support topology.
\endlemma

\proof
Let $\kappa$ be the kernel map of $\mathcal A$. Continuity of $\kappa$ with respect to the open support topology is, since $\kappa$ is continuous with respect to the lower Vietoris topology (because it is a kernel map), just the statement that for all $a\in A$ the set $\kappa^{-1}(\check a)$ is open. But, by Lemma~\ref{laxvsstrict}, this is equivalent to the open support property because
$\kappa^{-1}(\check a) = \{x\in X\st a\notin \kappa(x)\}$. 
\endproof

We shall denote by $\oSub A$ and $\oMax A$ the spaces $\Sub A$ and $\Max A$ equipped with the open support topology. 

\begin{corollary}
Let $X$ be a topological space whose soberification map $\sob:X\to\spectrum\topology(X)$ is surjective. Then the spectral vector bundles $(\pi,A,q)$ are, up to isomorphism, in bijective correspondence with the continuous maps $\kappa:X\to\oSub A$. This applies in particular to all sober spaces $X$, and thus also to Hausdorff spaces.
\end{corollary}

We note the following property, which will be used later on:

\lemma\label{oMaxT1}
$\oMax A$ is a $T_1$ space.
\endlemma

\proof
Let $V,W\in\oMax A$ be distinct elements. Without loss of generality assume that $V\not\subset W$, and let $a\in V\setminus W$. Then we have both $W\in\check a$ and $V\notin\check a$, on one hand, and, on the other, $W\notin\widetilde{A\setminus W}$ and $V\in\widetilde{A\setminus W}$. 
\endproof

Note, however, that $\oMax A$ is never a Hausdorff space if $\dim A\ge 2$. The following example illustrates this:

\example\label{oMaxnotHaus}
Each open set of $\oMax \CC^2$ is an open set of the lower Vietoris topology minus a finite number of rays. This implies that $\{(0,0)\}$ and $\linspan{(z,w)}$ cannot be separated by disjoint neighborhoods if $(z,w)\neq (0,0)$, since any basic open neighborhood of $\linspan{(z,w)}$ is generated by an open set of the Grassmannian $\Gr(1,\CC^2)$ (see \cite[Th.\ 7.5]{RS}).
\endexample

\paragraph{Fell topology.}

The open support topology is defined in a similar way to the Fell topology \cite{Fell62} (see also \cite{NT96,RS}), but the latter has subbasic open sets
\[
\check K=\{V\in\Max A\st V\cap K=\emptyset\}
\]
for each compact subset $K\subset A$ instead of only $\check a=\check{\{a\}}$ for $a\in A$. Therefore the Fell topology is obviously finer than the open support topology. Next we show that under mild restrictions it is strictly finer.

\lemma\label{fin2dim}
Let $A$ be a first countable topological vector space such that $\dim A>2$.
Then there is a compact set $K\subset A$
which is not contained in any finite union of 2-dimensional subspaces.
\endlemma

\proof
Let $e_0\in A$ and let $\{U_n\}$ be a countable basis of
neighborhoods of $e_0$ with $U_{n+1}\subset U_n$ for all $n$. There is
a sequence $(e_n)$ in $A$ such that $e_n\in U_n$ for all $n$, with the property that
for any $i<j<k$ the vectors $e_i,e_j,e_k$ are linearly independent:
assume we have chosen $e_1,\dots,e_{n-1}$ such that the list $e_0,\ldots,e_{n-1}$ has this property, and
let $C=\bigcup_{i,j<n}\linspan{e_i,e_j}$; since $A\setminus C$
is open and dense, we can choose $e_n\in U_n\cap(A\setminus C)$, and the list $e_0,\ldots,e_n$ has the required property.
Now let $K=\{e_n\}_{n\geq 0}$. The set $K$ is compact and
no 2-dimensional subspace contains more than two elements of $K$. 
\endproof

\theorem\label{fellfiner}
Let $A$ be a first countable topological vector space such that $\dim A>2$. Then the open support topology is strictly coarser
than the Fell topology.
\endtheorem

\proof
Let $K\subset A$ be a compact set which is not contained in any finite union of 2-dimensional subspaces (\cf\ Lemma~\ref{fin2dim}).
The set $\check K=\{V\in\Max A\st V\cap K=\emptyset\}$ is open in the Fell topology.
Let $a\in A$ be such that $\linspan a\in\check K$. We will show that no neighborhood of $\linspan a$ in the open support topology is contained in $\check K$.
Any basic neighborhood of $\linspan a$ in the open support topology is of the form
$\mathcal U\cap \check c_1\cap\dots\cap\check c_n$ with
$\mathcal U\subset\Max A$ open in the lower Vietoris topology and
$c_1,\dots,c_n\notin\linspan a$. But then $\linspan{a,b}\in\mathcal U$ for any $b\in A$. Let $C=\bigcup_i\linspan{a,c_i}$ and choose
$b\in K\setminus C$. Then for any $i$ we have
$b\notin\linspan{a,c_i}$, so $c_i\notin\linspan{a,b}$ because $c_i\notin\linspan a$.
It follows that
$\linspan{a,b}\in \mathcal U\cap \check c_1\cap\dots\cap\check c_n$.
But $\linspan{a,b}\notin\check K$, so $\check K$ is not open in the open support topology. 
\endproof

\subsection{Grassmannians}

\paragraph{The Grassmannian of a Hausdorff vector space.} Let $A$ be a Hausdorff vector space. 
Recall \cite{RS} that the \emph{$k$-Grassmannian} of $A$ is the set $\Gr(k,A)$ of $k$-dimensional linear subspaces of $A$, carrying the quotient topology given by the identification $\Gr(k,A)\cong\St(k,A)/\GL(k,\CC)$, where $\St(k,A)$ is the space of injective linear maps from $\CC^k$ to $A$ with the product topology.
Hence, the surjection $\St(k,A)\to\Gr(k,A)$ that is defined by $\phi\mapsto\image\phi$ is a continuous open map.

We have $\Gr(k,A)\subset\Max A$ because finite dimensional subspaces of a Hausdorff vector space are closed, and in \cite[Th.\ 7.2]{RS} it is shown that $\Gr(k,A)$ is a topological subspace of $\Max A$ with the lower Vietoris topology. In addition, below we shall use the fact \cite[Th.\ 7.5]{RS} that a basis of the lower Vietoris topology of $\Max A$ consists of all the sets of the form
\[
\upsegment\mathcal C := \{V\in\Max A\st \exists_{W\in\mathcal C}\ W\subset V\}\;,
\]
where $\mathcal C$ is an open set of a Grassmannian $\Gr(k,A)$ for some $k\in\NN_{>0}$.

\paragraph{Spaces of linear maps.}

Write $\LM(k,A)$ for the space of all the linear maps from $\CC^k$ to $A$ with the product topology. We note that we have a homeomorphism $\LM(k,A)\cong A^k$, under which $\St(k,A)$ can be identified with the subspace of linearly independent $k$-tuples.

\lemma\label{emptyint}
Let $A$ be a topological vector space. Any proper linear subspace $V\subset A$ has empty interior.
\endlemma
\proof
Let $V\subset A$ be a proper linear subspace. For all $U\subset V$, $v\in U$, and $w\notin V$, let $f\colon\CC\to A$
be the map defined by $f(\lambda)=v+\lambda w$. Then for all $\lambda\neq 0$ we have $f(\lambda)\notin V$, and thus $f^{-1}(U)=\{0\}$, showing that $U$ cannot be open. 
\endproof

\lemma\label{JPexla}
Let $A$ be a Hausdorff vector space. If $\dim A \geq k$ then $\St(k,A)$ is dense and open in $\LM(k,A)$.
\endlemma
\proof
We identify $\LM(k,A)$ with $A^k$.
Let us prove density. Let $U=U_1\times\dots\times U_k\subset A^k$ with each $U_i$ open in $A$, and let $\boldsymbol a=(a_1,\dots,a_k)\in U$ be such that $d=\dim\linspan{a_1,\dots,a_k}$ is maximal.
Assume that $d<k$. Then $\linspan{a_1,\dots,a_k}$ is a proper subspace of $A$, hence by Lemma~\ref{emptyint} with empty interior. Without loss of generality, assuming that $a_k\in\linspan{a_1,\ldots,a_{k-1}}$, we can choose $a_k'\in U_k$ such that
$\dim\linspan{a_1,\dots,a_{k-1},a_k'}>d$, a contradiction. Hence, $d=k$. This proves that $\St(k,A)$ intersects $U$, and thus it is dense in $\LM(k,A)$.
Now let $\boldsymbol a=(a_1,\dots,a_k)\in \St(k,A)$
be a linearly independent $k$-tuple. It has been shown in \cite[Lem.\ 7.1]{RS} that there are neighborhoods $U_i$ of $a_i$ such that every $(a_1',\dots,a_k')\in U_1\times\ldots\times U_k$ is linearly independent.
Such an open set $U_1\times\dots\times U_k$ is a neighborhood of $\boldsymbol a$
contained in $\St(k,A)$, and thus $\St(k,A)$ is open in $\LM(k,A)$. 
\endproof

\paragraph{Independent linear subspaces.}

We say that three linear subspaces $V_1,V_2,V_3\subset A$  are \emph{independent} if for all permutations $p$ of $\{1,2,3\}$ we have $V_{p(1)}\cap(V_{p(2)}+V_{p(3)})=\{0\}$. Equivalently, the mapping $(a,b,c)\mapsto a+b+c$ defines an isomorphism $V_1\oplus V_2\oplus V_3\cong V_1+ V_2+ V_3$.

\lemma\label{JPexl6}
Let $A$ be a Hausdorff vector space, $V\in\Max A$, and $k_1,k_2\in\NN_{>0}$.

\begin{enumerate}

\item\label{transvitem1} If $\codim V\geq k_1+k_2$ then the subspace of pairs $(W_1,W_2)$ such that
$V$, $W_1$ and $W_2$ are independent is open in $\Gr(k_1,A)\times\Gr(k_2,A)$.

\item\label{transvitem2} Let $a\notin V$. If $\codim V\geq k_1+k_2+1$ then the subset
\[\mathcal W_a=\{(W_1,W_2)\st a\notin V+W_1+W_2\quad\textrm{and}\quad V,W_1,W_2\textrm{ are independent}\}\]
is open and dense in $\Gr(k_1,A)\times\Gr(k_2,A)$.

\item\label{transvitem3} Let $\boldsymbol a=(a_1,\ldots,a_m)\in (A\setminus V)^m$. If $\codim V\geq k_1+k_2+1$ then the subset
\[\mathcal W_{\boldsymbol a}=\{(W_1,W_2)\st a_1,\ldots,a_m\notin V+W_1+W_2\quad\textrm{and}\quad V,W_1,W_2\textrm{ are independent}\}\]
is open and dense in $\Gr(k_1,A)\times\Gr(k_2,A)$.

\item\label{transvitem4}  Let $\boldsymbol a=(a_1,\ldots,a_m)\in (A\setminus V)^m$. If $\codim V\geq k_1+k_2+1$ then the subset
\[\mathcal X_{\boldsymbol a}=\{(W_1,W_2)\st a_1,\ldots,a_m\notin V+W_1+W_2\}\]
is dense in $\Gr(k_1,A)\times\Gr(k_2,A)$.

\end{enumerate}

\endlemma

\proof
Let $k=k_1+k_2$.
The quotient space $A/V$ is Hausdorff because $V$ is closed, and, since we have $\dim A/V\ge k$, by Lemma~\ref{JPexla} we conclude that $\St(k,A/V)$ is open in $\LM(k,A/V)$. Hence,
$\psi^{-1}(\St(k,A/V))$ is open in $\LM(k_1,A)\times\LM(k_2,A)$, where $\psi$ is given by the following composition:
\[
\xymatrix{
\LM(k_1,A)\times\LM(k_2,A)\ar@/^5ex/[rrrr]^{\psi}\ar[rr]_-{\cong}&&\LM(k,A)\ar[rr]_-{{((-)+V})_*}&&\LM(k,A/V)
}
\]
[Under the identification provided by the isomorphism $\LM(k,A)\cong A^k$, the map $\psi$ is the $k$-fold quotient map $A^k\to(A/V)^k$.]
Moreover, in fact we have $\psi^{-1}(\St(k,A/V))\subset\St(k_1,A)\times\St(k_2,A)$, and furthermore $\psi^{-1}(\St(k,A/V))$ is the subspace of pairs 
$(\phi_1,\phi_2)$ such that the spaces $V$, $\image\phi_1$ and $\image\phi_2$ are independent. Hence, since the quotient map $\St(k_i,A)\to\Gr(k_i,A)$ is open, the subspace of pairs $(W_1,W_2)\in \Gr(k_1,A)\times\Gr(k_2,A)$ such that $V$, $W_1$ and $W_2$ are independent is open in $\Gr(k_1,A)\times\Gr(k_2,A)$. This proves \eqref{transvitem1}.

In order to prove \eqref{transvitem2}, assume $a\notin V$ and $\codim V\ge k+1$, and let $V_a=V\oplus\linspan a$. Then $\dim A/V_a\ge k$ and, again by Lemma~\ref{JPexla}, $\St(k,A/V_a)$ is dense in $\LM(k,A/V_a)$. Since
$\psi$ is open and surjective, $\psi^{-1}\bigl(\St(k,A/V_a)\bigr)\subset\St(k_1,A)\times\St(k_2,A)$
is also dense. Now $\psi^{-1}\bigl(\St(k,A/V_a)\bigr)$ is the subspace of pairs $(\phi_1,\phi_2)\in\St(k_1,A)\times\St(k_2,A)$ such that the spaces $V_a$, $\image\phi_1$ and $\image\phi_2$ are independent. Such pairs are precisely the elements of the subspace
\[
\{(\phi_1,\phi_2)\st a\notin V+\image\phi_1+\image\phi_2\quad\textrm{and}\quad V,\image\phi_1,\image\phi_2\textrm{ are independent}\}\;,
\]
which therefore is open and dense in $\St(k_1,A)\times\St(k_2,A)$. Hence, applying the quotient $\phi_i\mapsto\image\phi_i$ we conclude that $\mathcal W_a$ is dense in $\Gr(k_1,A)\times\Gr(k_2,A)$, thus proving \eqref{transvitem2}.

Now let us prove \eqref{transvitem3}. For each $i$ we have $a_i\notin V$, by \eqref{transvitem2} the set $\mathcal W_{a_i}$ is open and dense in $\Gr(k_1,A)\times\Gr(k_2,A)$, and thus so is the finite intersection $\mathcal W_{\boldsymbol a} = \mathcal W_{a_1}\cap\ldots\cap\mathcal W_{a_m}$.

Finally, \eqref{transvitem4} follows from \eqref{transvitem3} because $\mathcal W_{\boldsymbol a}\subset\mathcal X_{\boldsymbol a}$.

\endproof

\subsection{Finer spectrum topologies}

\paragraph{Lower Vietoris topology revisited.}

In general it is not true, given some topological vector space $A$, that 
$\widetilde{U_1\cap U_2}=\widetilde U_1\cap\widetilde U_2$ for arbitrary open sets $U_1,U_2\subset A$. The following lemma provides an example where this equality holds, involving convex open sets.

\lemma\label{JPexl4}
Let $A$ be a topological vector space, let $U\subset A$ be a convex open set, and let $a\in A$. Let
\[
U^+=\bigcup_{\lambda\geq0}(U+\lambda a)\quad\text{and}\quad
U^-=\bigcup_{\lambda\geq0}(U-\lambda a)\,.
\]
Then $\widetilde U^+\cap\widetilde U^-=\widetilde U$.
\endlemma
\proof
We have $U^+\cap U^-=U$ so 
$\widetilde U\subset \widetilde U^+\cap\widetilde U^-$.
Let $V\in \widetilde U^+\cap\widetilde U^-$ and pick $c^+\in V\cap U^+$ and $c^-\in V\cap U^-$.
Then there are $b^-,b^+\in U$ and $\lambda^-,\lambda^+\geq 0$ such that 
$c^-=b^--\lambda^-a$ and $c^+=b^++\lambda^+a$. If $\lambda^+=0$ then $c^+=b^+\in V\cap U$ and thus $V\in\widetilde U$. If $\lambda^+\neq 0$
let $\mu=\lambda^-/(\lambda^++\lambda^-)\in[0,1]$. Then 
$\mu c^++(1-\mu)c^-=\mu b^++(1-\mu)b^-\in V\cap U$ since both $U$ and $V$ are convex,
so $V\in\widetilde U$. 
\endproof

\paragraph{Finer topologies and their prime open sets.}

Given $V_1,V_2\in\Max A$ with $V_1\subset V_2$, we introduce the following notation:
\[
\mathcal U_{V_1, V_2}=\Max A\setminus\bigl\{V\in\Max A\st V_1\subset V\subset V_2\bigr\}\,.
\]
Given a collection $\mathcal S$ of subsets of $\Max A$ we shall write $\SMax{\mathcal S}A$ for $\Max A$ equipped with the topology which is generated by the lower Vietoris topology and by $\mathcal S$.
We state the following lemma in more generality than needed in the rest of this paper. 

\lemma\label{JPexl0}
Let $A$ be a locally convex space, let $\mathcal S$ be a collection of subsets of $\Max A$, and
let $\mathcal Y\subset\SMax{\mathcal S} A$ be a subspace. Then for every prime open set $\mathcal P$
in $\mathcal Y$ there is $V\in\Max A$ such that for some union $S$ of elements of $\mathcal S$ we have
\[\mathcal P=\bigl(\mathcal U_{0,V}\cup S\bigr)\cap \mathcal Y\;.\]
\endlemma

\proof
Let $\mathcal P\in\spectrum\topology(\mathcal Y)$.
Since every prime open set is the union of subbasic open sets (\cf\ Lemma~\ref{JPexl2})
and for all families $(U_\alpha)$ of open sets of $A$ we have $\widetilde{\bigl(\bigcup_\alpha U_\alpha\bigr)}=\bigcup_\alpha\widetilde U_\alpha$,
there is $U\in\Omega(A)$ such that
$\mathcal P=(\widetilde U\cup S)\cap \mathcal Y$, 
where $S$ is a union of elements of $\mathcal S$.
Let $U_{\max}$ be the union of all the open sets $U_\alpha\in\Omega(A)$ such that 
$(\widetilde U_\alpha\cup S)\cap \mathcal Y=\mathcal P$. Then
\[
\bigl(\widetilde U_{\max}\cup S\bigr)\cap \mathcal Y
=\bigcup_\alpha\bigl(\widetilde U_\alpha\cup S\bigr)\cap \mathcal Y=\mathcal P\,.
\]
Let $V=A\setminus U_{\max}$. To complete the proof we only need to show that $V\in\Max A$, since then
$\widetilde U_{\max}=\widetilde{A
\setminus V}=\mathcal U_{0,V}$.

\begin{enumerate}

\item\label{step1} First we show that if $\lambda\in\CC$ and $a\in V$ then $\lambda a\in V$. 
The case $\lambda=0$ follows since $0\notin U_{\max}$, for otherwise we would have $\widetilde U_{\max}=\Max A$ and $\mathcal P=\mathcal Y$, which is a contradiction because $\mathcal P$ is prime.
So assume $\lambda\neq0$.
Then for any open set $U\in\Omega A$ we have $\widetilde{\lambda U}=\widetilde U$ and hence, 
by maximality, $\lambda U_{\max}=U_{\max}$. The result immediately follows.

\item Now we want to show that
$a,b\in V\Rightarrow a+b\in V$. We will prove this by contradiction. 
If $a+b\in U_{\max}$ then, by \eqref{step1}, we have $\frac12(a+b)\in U_{\max}$, so there is a convex open set $N\subset A$ satisfying $\frac12(a+b)\in N\subset  U_{\max}$.
Let
\[
N_a=\bigcup_{\lambda\geq0}\bigl(N+\lambda(a-b)\bigr)\quad\text{and}\quad N_b=\bigcup_{\lambda\geq0}\bigl(N+\lambda(b-a)\bigr)\,.
\]
Then $a\in N_a$ and $b\in N_b$, and, by Lemma~\ref{JPexl4}, we have 
\[
(\widetilde N_a\cap \mathcal Y)\cap(\widetilde N_b\cap \mathcal Y)=\widetilde N\cap \mathcal Y\subset\widetilde U_{\max}\cap \mathcal Y\subset\mathcal P\,.
\]
Hence, since $\mathcal P$ is prime, we must have either $\widetilde N_a\cap\mathcal Y\subset\mathcal P$ or $\widetilde N_b\cap\mathcal Y\subset\mathcal P$. Without loss of generality assume that $\widetilde N_a\cap \mathcal Y\subset\mathcal P$. Then
\[
\bigl((\widetilde{N_a\cup U_{\max}})\cup S\bigr)\cap \mathcal Y=\bigl(\widetilde N_a\cup\widetilde U_{\max}\cup S\bigr)\cap \mathcal Y=(\widetilde N_a\cap \mathcal Y)\cup\mathcal P=\mathcal P\;,
\]
and it follows that $N_a\subset U_{\max}$, and thus $a\in U_{\max}$. 
\end{enumerate}
\endproof

\subsection{Universal bundles}

\paragraph{Universal quotient vector bundles.}
We shall denote by
$U^1 A$
the universal bundle with Hausdorff fibers for $A$ (which is classified by the inclusion $\Max A\to\Sub A$). The following result shows that the open support property may fail even for bundles with Hausdorff fibers.

\theorem\label{JPext1}
Let $A$ be a locally convex space. The following conditions hold:
\begin{enumerate}
\item $\Max A$ is sober;
\item $U^1 A$ has continuous spectral kernel but it is not spectral because it does not satisfy the open support property.
\end{enumerate}
\endtheorem

\proof
Let $\mathcal Y=\Max A$ and let $\mathcal S=\emptyset$. Then Lemma~\ref{JPexl0} tells us that any prime open set is of the form
$\mathcal P=\mathcal U_{0,V}=\sob(V)$ with $V\in\Max A$, and it follows that the soberification map of $\Max A$ is surjective. It is also injective because $\Max A$ is $T_0$, and therefore $\Max A$ is sober. Now we claim that, for any $a\in A$,
we have $\osupp\hat a=\emptyset$.
Suppose we have $V\in\osupp\hat a$. Then $\hat a(V)\neq 0$, which is equivalent to $a\notin V$. Then, since open sets of $\Max A$ are upper closed we must also have $V\oplus \linspan a\in\osupp \hat a$, which implies $a\notin V\oplus \linspan a$, a contradiction. This proves the claim,
and it immediately follows that the universal bundle $U^1 A$ does not have the open support property. But its spectral kernel $\kfrak$ is continuous because it is constant, since for all $\mathcal P\in\Sigma\Omega(\Max A)$ we have
$\kfrak(\mathcal P)\supset\{a\in A\st \osupp\hat a\subset \mathcal P\}=A$. 
\endproof

\paragraph{Prime open sets for the open support topology.} Now we use the description of prime open sets of Lemma~\ref{JPexl0} in order to obtain a full characterization of the prime open sets for the open support topology.

\lemma\label{JPexl8}
Let $A$ be a topological vector space. Let also $V_1,V_2\in\Max A$ with $V_1\subset V_2$, let $\mathcal C\subset\Max A$, and let $a_1,\dots,a_k\in A$. The following conditions are equivalent:
\begin{enumerate}
\item $\upsegment\mathcal C\cap\check a_1\cap\dots\cap\check a_k\not\subset\mathcal U_{V_1,V_2}$;
\item There is $W\in \mathcal C$ such that $W\subset V_2$ and $a_1,\dots,a_k\notin W+V_1$.
\end{enumerate}
\endlemma
\proof
If $\upsegment\mathcal C\cap\check a_1\cap\dots\cap\check a_k\not\subset\mathcal U_{V_1,V_2}$ then there is
$V\in\Max A$ and $W\in\mathcal C$ such that $W\subset V$, $a_1,\dots,a_k\notin V$ and $V_1\subset V\subset V_2$.
It follows that $W\subset V_2$ and $a_1,\dots,a_k\notin W+V_1$. Conversely, if 
there is $W\in \mathcal C$ such that $W\subset V_2$ and $a_1,\dots,a_k\notin W+V_1$ then
$W+V_1\in \upsegment\mathcal C\cap\check a_1\cap\dots\cap\check a_k$ and $W+V_1\notin \mathcal U_{V_1,V_2}$. 
\endproof

\lemma\label{JPexl9.0}
Let $A$ be a Hausdorff vector space, and let $V_1,V_2\in\oMax A$ be such that
$V_1\subset V_2$ and the codimension of $V_1$ inside $V_2$ is infinite. Then $\mathcal U_{V_1,V_2}$ is a prime open set.
\endlemma

\proof
Let $\mathcal U_1$ and $\mathcal U_2$ be open sets of $\oMax A$. We shall prove that $\mathcal U_{V_1,V_2}$ is prime by showing that if $\mathcal U_1\not\subset\mathcal U_{V_1,V_2}$ and $\mathcal U_2\not\subset\mathcal U_{V_1,V_2}$
then $\mathcal U_1\cap\mathcal U_2\not\subset \mathcal U_{V_1,V_2}$. Assume then that
\begin{equation}\label{primehyp}
\mathcal U_1\not\subset\mathcal U_{V_1,V_2}\quad\textrm{and}\quad\mathcal U_2\not\subset\mathcal U_{V_1,V_2}\;.
\end{equation}
We may also assume that $\mathcal U_1$ and $\mathcal U_2$
are basic open sets, so, again by \cite[Th.\ 7.5]{RS}, we can write them in the form
\[
\mathcal U_1=\upsegment \mathcal C_1\cap \check a_1\cap\dots\cap \check a_k\quad \text{and}\quad
\mathcal U_2=\upsegment \mathcal C_2\cap \check b_1\cap\dots\cap \check b_n\;,
\]
where $\mathcal C_1\subset\Gr(d_1,A)$ and $\mathcal C_2\subset\Gr(d_2,A)$ are open sets of Grassmannians.
Then \eqref{primehyp} implies, by Lemma~\ref{JPexl8}, that $a_1,\ldots,a_k,b_1,\ldots,b_n\notin V_1$. Since the codimension of $V_1$ in $V_2$ is infinite, we can apply Lemma~\ref{JPexl6}\eqref{transvitem4} in order to conclude that the set
\[
\{(W_1,W_2)\in\Gr(d_1,V_2)\times\Gr(d_2,V_2)\st a_1,\ldots,a_k,b_1,\ldots,b_n\notin V_1+W_1+W_2\}
\]
is dense in $\Gr(d_1,V_2)\times\Gr(d_2,V_2)$ (notice that these are Grassmannians of $V_2$ rather than $A$). Condition \eqref{primehyp} implies that both $\mathcal U_1\cap\Gr(d_1,V_2)$ and $\mathcal U_2\cap \Gr(d_2,V_2)$ are not empty, so due to density there exist $W_1\in\mathcal U_1$ and $W_2\in\mathcal U_2$ such that $W_1+W_2\subset V_2$ and $a_1,\ldots,a_k,b_1,\ldots,b_n\notin V_1+W_1+W_2$.
Finally, again by Lemma~\ref{JPexl8}, we have
\[
\mathcal U_1\cap\mathcal U_2=\upsegment(\mathcal C_1+\mathcal C_2)\cap \check a_1\cap\dots\cap \check a_k\cap\check b_1\cap\dots\cap \check b_n\not\subset\mathcal U_{V_1,V_2}\;. 
\]
\endproof

\theorem\label{JPexl9}
Let $A$ be a Hausdorff locally convex space, and let $\mathcal P$ be a subset of $\oMax A$. Then $\mathcal P$ is a prime open set if and only if either
$\mathcal P=\sob(V)$ for some $V\in\oMax A$, or $\mathcal P=\mathcal U_{V_1,V_2}$ for some $V_1,V_2\in\oMax A$ such that
$V_1\subset V_2$ and the codimension of $V_1$ inside $V_2$ is infinite.
\endtheorem

\proof
Applying Lemma~\ref{JPexl0} with $\mathcal S=\{\check a\}_{a\in A}$ and $\mathcal Y=\Max^{\mathcal S} A=\oMax A$, we see that
$\mathcal P=\mathcal U_{0,V_2}\cup S$ with $V_2\in\oMax A$ and
$S$ a union of elements of $\mathcal S$. Then there is a subset $C\subset A$ such that
$S=\bigcup_{a\in C}\check a=\mathcal U_{V_1,A}$,
where $V_1=\overline{\spanmap C}$. Then
$
\mathcal P=\mathcal U_{0,V_2}\cup\mathcal U_{V_1,A}
$,
and, since $\mathcal P\neq \mathcal Y$ because $\mathcal P$ is prime, we must have $V_1\subset V_2$, and thus
$
\mathcal P=\mathcal U_{V_1,V_2}
$.
If $V_1= V_2$ then $\mathcal P=\mathcal Y\setminus\{V_1\}$, and this equals $\sob(V_1)$ because, by Lemma~\ref{oMaxT1}, $\oMax A$ is $T_1$. So let us assume $V_1\neq V_2$.
Suppose that the codimension $d$ of $V_1$ inside $V_2$ is finite.
If $d=1$ then $\mathcal U_{V_1,V_2}=\mathcal Y\setminus\{V_1,V_2\}$, which is not prime.
So let us assume that $d\ge 2$. Write $d=d_1+d_2$ with $d_1,d_2>0$,
and pick open sets $\mathcal C_1\subset\Gr(d_1,A)$ and $\mathcal C_2\subset\Gr(d_2,A)$ with the following properties (\cf\ Lemma~\ref{JPexl6}):
\begin{enumerate}
\item There are $W_1\in \mathcal C_1$ and $W_2\in\mathcal C_2$ such that $W_1+W_2\subset V_2$;
\item For any $W_1\in \mathcal C_1$ and $W_2\in\mathcal C_2$ the subspaces $W_1$, $W_2$ and $V_1$ are independent.
\end{enumerate}
Fix $W_1\in\mathcal C_1$ and $W_2\in\mathcal C_2$ with $W_1+W_2\subset V_2$ and pick $a\in V_2$ so that $a\notin V_1+ W_1$ and $a\notin V_1+W_2$. Due to \cite[Th.\ 7.5]{RS} the sets $\upsegment\mathcal C_1$ and $\upsegment\mathcal C_2$ are open in $\oMax A$, and from Lemma~\ref{JPexl8}
it follows that
\begin{equation}\label{antiprimeness}
\upsegment\mathcal C_1\cap\check a\not\subset\mathcal U_{V_1,V_2}\quad\textrm{and}\quad\upsegment\mathcal C_2\cap\check a\not\subset\mathcal U_{V_1,V_2}\;.
\end{equation}
Notice also that
\begin{equation}\label{antiprimeness2}
\bigl(\upsegment\mathcal C_1\cap\check a\bigr)\cap\bigl(\upsegment\mathcal C_2\cap\check a\bigr)=\upsegment(\mathcal C_1+\mathcal C_2)\cap\check a\;.
\end{equation}
By independence, for any $W\in \mathcal C_1+\mathcal C_2$ we have either $W\not\subset V_2$ or $W+ V_1=V_2$, and thus either $W\not\subset V_2$ or $a\in W+V_1$. Therefore, by Lemma~\ref{JPexl8}, we have
$\upsegment(\mathcal C_1+\mathcal C_2)\cap\check a\subset\mathcal U_{V_1,V_2}$, and, due to \eqref{antiprimeness} and \eqref{antiprimeness2}, we see that
$\mathcal U_{V_1,V_2}$ is not prime.
Finally, if the codimension of $V_1$ in $V_2$ is infinite then $\mathcal U_{V_1,V_2}$ is prime, by Lemma~\ref{JPexl9.0}. 
\endproof

\paragraph{Universal bundles versus open support topology.}

Let $A$ be a topological vector space. The identity $\ident:\oSub A\to\oSub A$ is the kernel map of a quotient vector bundle over $\oSub A$ that satisfies the open support property. However, as we shall see now, this is not in general a spectral vector bundle because the spectral kernel may fail to be continuous. In addition we show that the open support property does not in general entail having a closed zero section.

We shall denote by $U^\circ A$ the bundle which is classified by the inclusion $\oMax A\to\oSub A$.

\theorem\label{JPext2}
Let $A$ be a Hausdorff vector space.
\begin{enumerate}
\item If $\dim A=\infty$ the open support property for $U^{\circ} A$ holds but its spectral kernel is not continuous.
\item If $\dim A<\infty$ then $U^{\circ} A$ is sober (and thus spectral).
\end{enumerate}
\endtheorem

\proof
Let $\kappa$ and $\kfrak$ be the kernel map and the spectral kernel map of $U^{\circ} A$ ($\kappa$ is the identity on $\Max A$). Then $\kappa=\kfrak\circ\sob$.
Assuming first that $\dim A=\infty$,
from Lemma~\ref{JPexl9.0} it follows that $\emptyset=\mathcal U_{0,A}$ is a prime open set of $\oMax A$. Given any $V\in\oMax A$ we have
$\kfrak(\sob(V))=\kappa(V)=V$.
If $\kfrak$ were continuous then from Theorem~\ref{thmcontkfrak} we would conclude, since $\emptyset\subset V$ for all $V\in\oMax A$, that $\kfrak(\emptyset)=V$ for all $V\in\oMax A$, a contradiction.
If $A$ has finite dimension then it is Euclidean and, by Theorem~\ref{JPexl9},
the only prime open sets are of the form $\sob(V)$ for $V\in\oMax A$. This means that $\oMax A$ is sober, and, by Lemma~\ref{surjsob}, we conclude that $\kfrak$ is continuous. Hence, for finite dimensional $A$ the bundle $\mathcal U^{\circ} A$ is sober.

\endproof

In addition, we mention the following simple consequence of Theorem~\ref{fellfiner}:

\begin{corollary}
Let $A$ be a first countable topological vector space such that $\dim A>2$. Then the zero section of $U^{\circ} A$ is not closed.
\end{corollary}

\proof
By \cite[Th.\ 5.7]{RS}, the kernel of a quotient vector bundle with closed zero section must be Fell continuous, so the result follows immediately from Theorem~\ref{fellfiner}. 
\endproof


\begin{thebibliography}{10}

\bibitem{BCL}
Paolo Bertozzini, Roberto Conti, and Wicharn Lewkeeratiyutkul.
\newblock A horizontal categorification of {G}el'fand duality.
\newblock {\em Adv. Math.}, 226(1):584--607, 2011.

\bibitem{Bredon}
Glen~E. Bredon.
\newblock {\em Sheaf theory}, volume 170 of {\em Graduate Texts in
  Mathematics}.
\newblock Springer-Verlag, New York, second edition, 1997.

\bibitem{BE12}
Alcides Buss and Ruy Exel.
\newblock Fell bundles over inverse semigroups and twisted \'etale groupoids.
\newblock {\em J. Operator Theory}, 67(1):153--205, 2012.

\bibitem{BM16}
Alcides Buss and Ralph Meyer.
\newblock Inverse semigroup actions on groupoids.
\newblock {\em Rocky Mountain J. Math.}, 2016.
\newblock To appear.

\bibitem{Exel}
Ruy Exel.
\newblock Noncommutative {C}artan subalgebras of {$C^*$}-algebras.
\newblock {\em New York J. Math.}, 17:331--382, 2011.

\bibitem{FeldmanMooreI-II}
Jacob Feldman and Calvin~C. Moore.
\newblock Ergodic equivalence relations, cohomology, and von {N}eumann
  algebras. {I, II}.
\newblock {\em Trans. Amer. Math. Soc.}, 234(2):289--359, 1977.

\bibitem{Fell62}
J.~M.~G. Fell.
\newblock A {H}ausdorff topology for the closed subsets of a locally compact
  non-{H}ausdorff space.
\newblock {\em Proc. Amer. Math. Soc.}, 13:472--476, 1962.

\bibitem{FD1}
J.~M.~G. Fell and R.~S. Doran.
\newblock {\em Representations of {$^*$}-algebras, locally compact groups, and
  {B}anach {$^*$}-algebraic bundles. {V}ol. 1}, volume 125 of {\em Pure and
  Applied Mathematics}.
\newblock Academic Press, Inc., Boston, MA, 1988.
\newblock Basic representation theory of groups and algebras.

\bibitem{stonespaces}
Peter~T. Johnstone.
\newblock {\em Stone spaces}, volume~3 of {\em Cambridge Studies in Advanced
  Mathematics}.
\newblock Cambridge University Press, Cambridge, 1986.
\newblock Reprint of the 1982 edition.

\bibitem{JT}
Andr{\'e} Joyal and Myles Tierney.
\newblock An extension of the {G}alois theory of {G}rothendieck.
\newblock {\em Mem. Amer. Math. Soc.}, 51(309):vii+71, 1984.

\bibitem{KR}
David Kruml and Pedro Resende.
\newblock On quantales that classify {$C^\ast$}-algebras.
\newblock {\em Cah. Topol. G\'eom. Diff\'er. Cat\'eg.}, 45(4):287--296, 2004.

\bibitem{Kumjian98}
Alex Kumjian.
\newblock Fell bundles over groupoids.
\newblock {\em Proc. Amer. Math. Soc.}, 126(4):1115--1125, 1998.

\bibitem{Kumjian}
Alexander Kumjian.
\newblock On {$C^\ast$}-diagonals.
\newblock {\em Canad. J. Math.}, 38(4):969--1008, 1986.

\bibitem{MP1}
Christopher~J. Mulvey and Joan~Wick Pelletier.
\newblock On the quantisation of points.
\newblock {\em J. Pure Appl. Algebra}, 159(2-3):231--295, 2001.

\bibitem{MP2}
Christopher~J. Mulvey and Joan~Wick Pelletier.
\newblock On the quantisation of spaces.
\newblock {\em J. Pure Appl. Algebra}, 175(1-3):289--325, 2002.
\newblock Special volume celebrating the 70th birthday of Professor Max Kelly.

\bibitem{NT96}
Tsugunori Nogura and Dmitri Shakhmatov.
\newblock When does the {F}ell topology on a hyperspace of closed sets coincide
  with the meet of the upper {K}uratowski and the lower {V}ietoris topologies?
\newblock In {\em Proceedings of the {I}nternational {C}onference on
  {C}onvergence {T}heory ({D}ijon, 1994)}, volume~70, pages 213--243, 1996.

\bibitem{Paterson}
Alan L.~T. Paterson.
\newblock {\em Groupoids, inverse semigroups, and their operator algebras},
  volume 170 of {\em Progress in Mathematics}.
\newblock Birkh\"auser Boston, Inc., Boston, MA, 1999.

\bibitem{PR12}
M.~Clarence Protin and Pedro Resende.
\newblock Quantales of open groupoids.
\newblock {\em J. Noncommut. Geom.}, 6(2):199--247, 2012.

\bibitem{RenaultLNMath}
Jean Renault.
\newblock {\em A groupoid approach to {$C^{\ast} $}-algebras}, volume 793 of
  {\em Lecture Notes in Mathematics}.
\newblock Springer, Berlin, 1980.

\bibitem{Renault}
Jean Renault.
\newblock Cartan subalgebras in {$C^*$}-algebras.
\newblock {\em Irish Math. Soc. Bull.}, 61:29--63, 2008.

\bibitem{Re07}
Pedro Resende.
\newblock \'{E}tale groupoids and their quantales.
\newblock {\em Adv. Math.}, 208(1):147--209, 2007.

\bibitem{RS}
Pedro Resende and Jo\~{a}o Santos.
\newblock Open quotients of trivial vector bundles.
\newblock arXiv:1510.06329.

\bibitem{RV}
Pedro Resende and Steven Vickers.
\newblock Localic sup-lattices and tropological systems.
\newblock {\em Theoret. Comput. Sci.}, 305(1-3):311--346, 2003.
\newblock Topology in computer science (Schlo{\ss} Dagstuhl, 2000).

\bibitem{Sieben}
N.~Sieben.
\newblock Fell bundles over $r$-discrete groupoids and inverse semigroups.
\newblock Unpublished draft, available at http://jan.ucc.nau.edu/$\sim
  $ns46/bundle.ps.gz.

\bibitem{Steenrod}
Norman Steenrod.
\newblock {\em The topology of fibre bundles}.
\newblock Princeton Landmarks in Mathematics. Princeton University Press,
  Princeton, NJ, 1999.
\newblock Reprint of the 1957 edition, Princeton Paperbacks.

\bibitem{Vietoris}
Leopold Vietoris.
\newblock Bereiche zweiter {O}rdnung.
\newblock {\em Monatsh. Math. Phys.}, 32:258--280, 1922.

\end{thebibliography}
\end{document}